\documentclass[12 pt,reqno]{amsart}

\usepackage{latexsym, amsmath, amssymb, longtable, booktabs,amscd,microtype,booktabs,cases}
\usepackage[square,numbers]{natbib}
\usepackage{graphicx}
\usepackage[dvipsnames]{xcolor}
\usepackage{caption}
\usepackage{dsfont}
\usepackage[all]{xy}
\usepackage{enumerate}
\usepackage{multirow}
\usepackage{tikz}
\usetikzlibrary{matrix,arrows,decorations.pathmorphing}
\usepackage{collectbox}
\usepackage{enumerate}
\usepackage[colorinlistoftodos,prependcaption, textsize=tiny]{todonotes}
\reversemarginpar
\usepackage{amsmath}
\usepackage{enumitem}
\usepackage{hyperref}
\let\oldtocsection=\tocsection

\let\oldtocsubsection=\tocsubsection

\let\oldtocsubsubsection=\tocsubsubsection

\renewcommand{\tocsection}[2]{\hspace{0em}{\vspace{0.5em}}\oldtocsection{#1}{#2}}
\renewcommand{\tocsubsection}[2]{\hspace{1em}{\vspace{0.5em}}\oldtocsubsection{#1}{#2}}
\renewcommand{\tocsubsubsection}[2]{\hspace{2em}\oldtocsubsubsection{#1}{#2}}
\hypersetup{
	colorlinks,
	citecolor=red,
	filecolor=black,
	linkcolor=blue,
	urlcolor=black
}

\numberwithin{equation}{section}
\newcommand{\bigzero}{\mbox{\normalfont\Large\bfseries 0}}
\newcommand{\bigi}{\mbox{\normalfont\Large\bfseries I}}
\newcommand{\rvline}{\hspace*{-\arraycolsep}\vline\hspace*{-\arraycolsep}}

\newcommand{\Z}{\mathbb{Z}}

\newcommand\fh{\mathfrak{h}}
\newcommand\CT{\mathcal{T}}

\newcommand{\C}{\mathbb{C}}

\newcommand\N{\mathbb{N}}

\newcommand{\m}{\textbf{m}}
\newcommand{\um}{\overline{\textbf{m}}}
\newcommand{\e}{\mathfrak{a}}
\newcommand{\bb}{\mathfrak{b}}
\newcommand{\V}{\mathcal{V}}
\newcommand{\lk}{\underbar{\textbf{k}}}
\newcommand{\lm}{\underbar{\textbf{m}}}
\newcommand{\uk}{\overline{\textbf{k}}}
\newcommand{\ukk}{\overline{\textbf{k}^{\prime}}}

\makeatother

\newtheorem{thm}{Theorem}[section]
\newtheorem{theorem}[thm]{Theorem}
\newtheorem{cor}[thm]{Corollary}

\newtheorem{prop}[thm]{Proposition}
\newtheorem{lemma}[thm]{Lemma}
\theoremstyle{definition}
\newtheorem{definition}[thm]{Definition}
\newtheorem{remark}[thm]{Remark}

\theoremstyle{definition}

\theoremstyle{remark}

\theoremstyle{remark}

\makeatletter
\def\imod#1{\allowbreak\mkern10mu({\operator@font mod}\,\,#1)}
\makeatother
\begin{document}
	
\title{Weyl modules for toroidal Lie algebras}
\author{ Sudipta Mukherjee, Santosha Kumar Pattanayak, Sachin S. Sharma}
\date{}

\maketitle
\begin{abstract}
In this paper we study Weyl modules for a toroidal Lie algebra $\CT$ with arbitrary $n$ variables. Using the work of Rao \cite{1995}, we prove that the level one global Weyl modules of $\CT$
are isomorphic to suitable submodules of a Fock space representation of $\CT$ upto a twist. As an application, we compute the graded character of the level one 
local Weyl module of $\CT$, thereby generalising the work of Kodera \cite{ko}.

\end{abstract}
\tableofcontents
\section{Introduction}
The notion of Weyl modules emanated from the work of Chari and Pressley \cite{CP}. Since then this topic is actively pursued by many mathematicians \cite{CK, cv, 2002, CL, SS, FL, khandai, Na, san}. Weyl modules (global and local) were initially defined for loop algebras. 
Global Weyl modules are the maximal integrable highest weight modules for the loop algebras and their maximal finite dimensional integrable quotients are local Weyl modules. The graded characters of local Weyl modules gained significant attention when Chari and Loktev \cite{CL} proved that the local Weyl modules,  Demazure modules and fusion product modules over the current algebra $\mathfrak{sl}_n[t]$ are all isomorphic.
So the graded characters of local Weyl modules can also be expressed as Kostka Foulkes polynomials. Chari-Le \cite{2002} made the first attempt of generalising these notions to the quotient of toroidal Lie algebras in two variables and in 
\cite{SS} it was generalised
for loop Kac-Moody algebras.Very recently, Kodera \cite{ko} defined Weyl modules for toroidal Lie algebras in two variables and building on work of \cite{yokunama}, he obtained realisations of  level one global Weyl modules with modules constructed in \cite{yokunama}. He further obtained the graded character of level one local Weyl module. In this paper we generalise Kodera's  work for toroidal Lie algebras in arbitrary $n$ variables. More precisely, we extend the notion of global and local Weyl modules for toroidal Lie algebras in arbitrary $n$ variables. We use the work of Rao \cite{1995} in which sub-quotients of  Fock space representation of these toroidal Lie algebras are studied using vertex operators.  

\subsection{Presentation} In section \ref{s1} we give the definition of a toroidal Lie algebra $\CT$ and its realisation in terms of generators and relations as given in \cite{CMP}. We define an important subalgebra $\CT^{+}$ of $\CT$ which play a key role in this work. We end this section with triangular decomposition of $\CT$. We start Section \ref{s2} with definitions of global and local Weyl module for $\CT$. For a given dominant integral weight $\Lambda$, one associates the $\CT$ module $W_{\mathrm{glob}}{(\Lambda)}$. The $\CT$-module $W_{\mathrm{glob}}{(\Lambda})$ becomes 
right module for commutative algebra $A{(\Lambda)}$, the algebra associated with $W_{\mathrm{glob}}{(\Lambda)}$.  Using Garland identities and results in  \cite{CK} we prove that $A(\Lambda)$ is isomorphic to a ring of invariants of a subgroup of the symmetric group (Theorem \ref{t1}). The local Weyl modules $W_{\mathrm{loc}}(\Lambda)$ are tensor products of global Weyl modules with one dimensional representations of $A(\Lambda)$. All the analogous definitions are defined similarly for $\CT^+$. We end Section \ref{s2} with a result which implies that the local Weyl modules for $\CT$ are generated by $\mathcal{U}(\CT^+).v_{\Lambda}$, where $v_{\Lambda}$ is a highest weight vector of $W_{\mathrm{loc}}$.
In section \ref{s3}, we utilise results in Section \ref{s2} to the special case $\Lambda = \Lambda_0$ to obtain an upper bound for the local Weyl module of $\CT^+$ and hence $\CT$. Finally, in Section \ref{s4} we begin with recalling basics of Fock space representation of $\CT$ using vertex operators.
We invoke results in \cite{1995} in which, Rao studied the sub-quotients of Fock space representation of $\CT$ and determined $\CT$ action on them using vertex operators, to prove that the level one global Weyl modules are isomorphic to suitable submodules of  the Fock space of $\CT$ up to a twist. 

\vspace{0.5cm}

Let $\CT$ be a toroidal Lie algebra in $n$ variables.  Then the following are our main results (See Section \ref{s4} for details).
\begin{thm} As a $\CT$- module $W_{\mathrm{glob}}(\Lambda_0)$ is isomorphic to $A^* V(0)$.
\end{thm} 
\begin{thm}
The graded character of $W_{\mathrm{loc}}(\Lambda_0)$ is given by
	$$\mathrm{ch}_{q_1, q_2, \ldots, q_{n}} W_{\mathrm{loc}}(\Lambda_0) = \hspace{0.2cm} {\displaystyle{ \mathrm{ch}_{q_1} L(\Lambda_0)}\prod_{m>0, i = 2}^{n}{\frac{1}{1- {q_1}^m q_i}}}.$$
\end{thm}

\section{Notation}
Throughout the paper $\C$, $\Z$, $\N$, $\Z_{+}$, $\Z_{-}$  denote the set of complex numbers, integers, natural numbers , non-negative integers and non-positive integers respectively. Let $I$ (respectively $I_{\mathrm{aff}})$ denote $\{1,2,\cdots,l\}$ (respectively $\{0,1,2,\cdots,l\})$. We denote $\textbf{m}$, $\underbar{\textbf{m}} $, $\um$  by $(m_1,\cdots,m_n) \in \Z^n$, $(m_2,\cdots,m_n) \in \Z^{n-1}$ and $(m_1,\cdots,m_{n-1}) \in \Z^{n-1}$ respectively. We denote the universal enveloping algebra of a Lie algebra $\mathfrak g$ by $\mathcal{U}(\mathfrak g)$.

	\section{Preliminaries} \label{s1}
\subsection{Basics of toroidal Lie algebras}
In this section we recall the definition and the explicit realization of toroidal Lie algebras. Throughout this paper we work with the base field  $\mathbb{C}$.
Let $\mathfrak{g}$ be a finite dimensional simple Lie algebra over $\mathbb C$ of rank $l$ with Cartan subalgebra $\mathfrak{h}$. Let $R_{fin}$ be the set of roots of $\mathfrak{g}$ with respect to $\mathfrak{h}$. Fix a set of simple roots $\alpha_{i}$ ($i \in I$), simple coroots $\alpha_{i}^{\lor}$ ($i \in I)$ and  fundamental weights $\Lambda_{i}$ ($i \in I)$ of $\mathfrak{g}$ with respect to $\mathfrak{h}$. Let $R_{fin}^{+}$ be the set of positive roots  with respect to a fixed Borel subalgebra containing $\mathfrak h$ and let $\theta \in R_{fin}^{+}$ be the highest root of $\mathfrak{g}$. Let $Q_{fin}=\oplus_{i=1}^l \Z \alpha_i$ be the root lattice of $\mathfrak{g}$. For a given $\alpha \in R_{fin}^{+}$, let $\mathfrak{g}_{\pm \alpha}$ be the corresponding root spaces and fix non-zero elements $x_{\pm \alpha} \in \mathfrak{g}_{\pm \alpha}$, $h_{\alpha} \in \mathfrak{h}$ such that
\begin{center}
$[x_{\alpha},x_{-\alpha}]=h_{\alpha}$,  \hspace{0.6cm}   $[h_{\alpha},x_{\pm\alpha}]=2x_{\pm\alpha}$.
\end{center}\vspace{0.2cm}
Set $x_{\alpha_{i}}=e_{i}$, $x_{-\alpha_{i}}=f_{i}$, $h_{\alpha_{i}}=h_{i}$ so that the subalgebra generated by $\{e_{i},f_{i},h_{i}\}$ is isomorphic to $\text{sl}_{2}(\mathbb C)$. Let $(,)$ be a non-degenerate symmetric bilinear form on $\mathfrak{g}$ normalized by $(\theta, \theta)=2$. 
\par
\vspace{0.2cm}
For a positive integer $n$, let $A_{n}=\C[t_{1}^{\pm 1} ,\cdots, t_{n}^{\pm 1}]$ be the Laurent polynomial ring in $n$ variables $t_1, t_2, \cdots, t_n$. For $\textbf{m} \in \Z^n$ (resp. $\underbar{\textbf{m}} \in \Z^{n-1}$, $\um \in \Z^n$), let $t^{\textbf{m}}$ (resp. $t^{\underbar{\textbf{m}}}$, $t^{\um}$) denote the element $t_1^{m_1} \cdots t_n^{m_n}$(resp. $t_2^{m_2} \cdots t_n^{m_n}$ , $t_1^{m_1} \cdots t_{n-1}^{m_{n-1}}$) in $A_{n}$.
Let $L(\mathfrak{g})= \mathfrak{g}\otimes \C[t_{1}^{\pm 1} ,\cdots, t_{n}^{\pm 1}]$ and $\mathcal{Z}=\Omega/dL$ be the space of K\"ahler differentials spanned by the vectors $\{t^{\textbf{m}}K_{i}: \textbf{m} \in \Z^n, 1\leq i \leq n\}$ subject to the relation $\sum_{i=1}^n m_{i}\,t^{\textbf{m}}K_{i}=0$. For $1\leq i \leq n$, we have $n$ derivation on
 $L(\mathfrak{g})\oplus \mathcal{Z}$ given by
 \vspace{0.1cm}
 
 \begin{center}
 	$d_{i}(x \otimes t^{\textbf{m}})=m_{i}x \otimes t^{\textbf{m}}$,\,\,\,  $d_{i}(t^{\textbf{m}}K_{j})=m_{i}t^{\textbf{m}}K_{j}$, \,\,\,\,\, $x\in \mathfrak{g}$, $\textbf{m} \in \Z^{n}$, $1\leq i,j \leq n$.
 \end{center} 
\vspace{0.2cm}

Let $D$ be the $\C$ linear span of $d_1,d_2,\cdots,d_n$. Then the algebra $\mathcal{T}=L(\mathfrak{g})\oplus \mathcal{Z} \oplus D $ is called a toroidal Lie algebra. So, as a vector space 
$$\mathcal{T}=\mathfrak{g} \otimes \C[t_1 ^{\pm},t_2^{\pm},\cdots, t_n^{\pm}] \hspace{0.2cm} \oplus \sum\limits_{\substack{i=1 \\ \textbf{m} \in \Z^{n}}}^n \C t^{\textbf{m}} K_i  \oplus \hspace{0.1cm} \bigoplus\limits_{i=1}^n \C d_i,$$
and the Lie bracket on $\mathcal{T}$ is given by:
\begin{center}
	$[x \otimes t^{\textbf{p}},\hspace{0.1cm} y \otimes t^{\textbf{q}}]=[x,y] \otimes t^{\textbf{p+q}}+(x,y)\hspace{0.05cm} \displaystyle\sum_{i=1}^n p_{i}\,t^{\textbf{p
			+q}}K_{i},$
\end{center}
 \begin{center}
	 $[x \otimes t^{\textbf{p}},z]=0$, \hspace{0.2cm} $[z,z']=0$, \hspace{0.2cm} $[d_i,d_j]=0,$
\end{center}
\begin{center}
	$[d_{i},x \otimes t^{\textbf{p}}]=p_{i}x \otimes t^{\textbf{p}}$, \hspace{0.2cm} $[d_{i},t^{\textbf{p}}K_{j}]=p_{i}t^{\textbf{p}}K_{j},$
\end{center}
where $x,y \in \mathfrak{g}$, $\textbf{p}, \textbf{q} \in \Z^{n}$, $z,z' \in \mathcal{Z}$, $1\leq i,j \leq n$  and $(,)$ is the non-degenerate symmetric bilinear form on $\mathfrak{g}$ defined above. Let $\mathcal{Z}_0$ be the subspace of $\mathcal{Z}$ spanned by the zero degree central elements $\{K_i:1\leq i \leq n\}$ and let $\mathfrak{h}_{\mathcal{T}}=\mathfrak{h}\oplus \mathcal{Z}_0 \oplus D$. We identify $\mathfrak{h}^{*}$ as a subspace of $\mathfrak{h}_{\mathcal{T}}^{*}$ and an element $\lambda \in \mathfrak{h}^{*}$ is extended to an element of $\mathfrak{h}_{\mathcal{T}}^{*}$ by setting $\lambda(K_i)=\lambda(d_i)=0$ for $1\leq i\leq n$. For $1\leq i\leq n$, define $\delta_{i} \in \mathfrak{h}_{\mathcal{T}}^{*}$ by 
\vspace{0.1cm}
\begin{center}
	$\delta_{i}\,({\mathfrak{h}\oplus \mathcal{Z}_0 })=0$, \hspace{0.5cm}  $\delta_{i}(d_j)=\delta_{ij}$ for $1\leq j\leq n$.
\end{center}
\vspace{0.1cm}
 For $\alpha \in R_{fin}$ and $\textbf{m}=(m_1,\cdots,m_n) \in \Z^n$, set $\delta_{\textbf{m}}=\sum_{i=1}^n m_{i}\, \delta_{i}$. Let $R_{\mathcal{T}}$ (resp. $R_{\mathcal{T}}^{re}$ and $R_{\mathcal{T}}^{im}$ ) be the set of all (resp. real and imaginary) roots of $\mathcal{T}$ with respect to $\mathfrak{h}_{\mathcal{T}}$. Then we have
 \begin{center}
 	$R_{\mathcal{T}}^{re}=\{\alpha+\delta_{\textbf{m}}: \alpha \in R_{fin}, \textbf{m} \in \Z^{n}\}$,  \hspace{0.2cm}	$R_{\mathcal{T}}^{im}=\{\delta_{\textbf{m}}: \textbf{m} \in \Z^{n}\setminus \{\textbf{0}\}\}$ and $R_{\mathcal{T}}=R_{\mathcal{T}}^{re} \cup R_{\mathcal{T}}^{im}.$
 \end{center} The root spaces for $\mathcal{T}$ are given by
\begin{center}
	$\mathcal{T}_{\alpha+\delta_{\textbf{m}}}=\mathfrak{g}_{\alpha} \otimes t^{\textbf{m}}$, $\alpha \in R_{fin}$,\\
	\vspace{0.1cm}
	$\mathcal{T}_{\delta_{\textbf{m}}}=\mathfrak{h} \otimes t^{\textbf{m}} \oplus \displaystyle\sum_{i=1}^n \mathbb{C} t^{\textbf{m}}K_{i}$, \hspace{0.3cm}    $\textbf{m}\neq \textbf{0},$\\
	\vspace{0.1cm}
	$\mathcal{T}_{0}=\mathfrak{h}_{\mathcal{T}}.$
	\end{center}
Set $\alpha_{l+i}:=\delta_{i}-\theta$ for $ 1 \leq i \leq n$, then $\{\alpha_1, \cdots, \alpha_l, \alpha_{l+1},\cdots, \alpha_{l+n}\}$ forms a simple root system for $R_{\mathcal{T}}$. For a real root $\beta=\alpha+\delta_{\textbf{m}}$ with $\alpha \in R_{fin}^{+}$ and $\textbf{m} \in \Z^{n}$, set $\beta^{\lor}=\alpha^{\lor}+\frac{2}{(\alpha, \alpha)} \sum_{i=1}^n m_{i}K_{i}$. Then with the Lie bracket on $\mathcal{T}$, we see that the subalgebra spanned by $\{x_{\alpha} \otimes t^{\textbf{m}},\, x_{- \alpha} \otimes t^{-{\textbf{m}}},\, \beta^{\lor} \}$ is isomorphic to $\mathfrak{sl}_2(\C)$.\par
\vspace{0.2cm}
The toroidal Lie algebra $\mathcal{T}$ contains a Lie subalgebra $\mathcal{T}_{\mathrm{aff}}$ isomorphic to the untwisted affine Lie algebra associated with $\mathfrak{g}$:
\begin{center}
	$\mathcal{T}_{\mathrm{aff}}=\mathfrak{g}\otimes \C[t_{1}^{\pm 1}] \oplus \C K_1 \oplus \C d_1$.
\end{center}
\vspace{0.1cm}
A Cartan subalgebra of $\mathcal{T}_{\mathrm{aff}}$ is given by $\mathfrak{h}_{\mathrm{aff}}=\mathfrak{h} \, \oplus \, \C K_1 \oplus \, \C d_1$. With respect to $\mathfrak{h}_{\mathrm{aff}}$, the roots of $\mathcal{T}_{\mathrm{aff}}$ are described by the set $R_{\mathrm{aff}}= \{\alpha+ m \delta_1: \alpha \in R_{fin}, \, m \in \Z\} \cup \{m\delta_1: m \in \Z \setminus \{0\}\}$.  
The real and imaginary roots of $\mathcal{T}_{\mathrm{aff}}$ are described by the first and second set on the expression of $R_{\mathrm{aff}}$ respectively.
Let $R_{\mathrm{aff}}^{+}$ and $R_{\mathrm{aff}}^{-}$ denote the set of positive and negative roots respectively. We have
\begin{center}
	$R_{\mathrm{aff}}^{\pm}=\{\alpha+m\delta_{1}: \alpha \in R_{fin}, m \in \Z_{\pm} \setminus \{0\}\}\cup \{m\delta_1: m \in \Z_{\pm} \setminus \{0\}\} \cup  R_{fin}^{\pm}$.
\end{center}
\vspace{0.1cm}
 The set of simple roots \,$\Delta_{\mathrm{aff}}$ \,and the simple coroots \,$\Delta_{\mathrm{aff}}^{\lor}$\, of \,$\mathcal{T}_{\mathrm{aff}}$ are given by $\Delta_{\mathrm{aff}}=\{\alpha_1,\cdots,\alpha_l,\alpha_{l+1}=\delta_{1}-\theta\}$ and $\Delta_{\mathrm{aff}}^{\lor}=\{\alpha_1^{\lor},\cdots,\alpha_l^{\lor},\alpha_{l+1}^{\lor}=K_{1}-\theta^{\lor}\}$ respectively. For a root $\beta \in R_{\mathrm{aff}}$, let 
$\mathcal{T}_{\mathrm{aff}}^{\beta}$ denote the corresponding root space of $\mathcal{T}_{\mathrm{aff}}$. We set $n_{\mathrm{aff}}^{\pm}=\displaystyle\oplus_{\beta \in R_{\mathrm{aff}}^{\pm}} \mathcal{T}_{\mathrm{aff}}^{\beta}$. From now on we denote $\alpha_{l+1}$, $\alpha_{l+1}^{\lor}$ by $\alpha_0$ and $\alpha_{0}^{\lor}$ respectively. Let $\{\Lambda_{0}, \Lambda_{1},\cdots, \Lambda_{l}\}$ be the set of fundamental weights of $\mathcal{T}_{\mathrm{aff}}$, that is, $\Lambda_{i}(\alpha_{j}^{\lor})= \delta_{ij}$ for $i,j \in I_{\mathrm{aff}}$. Let $L(\lambda)$ denote the irreducible highest weight module generated by the highest weight vector $v_{\lambda}$ with highest weight $\lambda$. Defining relations in $L(\lambda)$ are given by
\begin{center}
	$n_{\mathrm{aff}}^{+}.\, v_{\lambda}=0$, \hspace{0.6cm} $h. v_{\lambda}= \lambda(h)  v_{\lambda}$,\,\,\, $\forall h \in \mathfrak{h}_{\mathrm{aff}}$, \hspace{0.6cm} $ f_i^{\lambda(h_i)+1} . v_{\lambda}=0$\,\,\, for $i \in I_{\mathrm{aff}}$.
\end{center}\par
\vspace{0.3cm}
Now we recall some identities by Garland in the universal enveloping algebra of $\mathcal{T}_{\mathrm{aff}}$. For $\alpha \in R_{fin}$, define a power series $\textbf{P}_{\alpha}(u)$ in indeterminate u with coefficients in $\mathcal{U}(\mathfrak{h} \otimes \C[t_1^\pm])$:
\vspace{0.2cm}
$$\textbf{P}_{\alpha}(u)=\mathrm{exp} \hspace{0.1cm}\Bigg(-\displaystyle\sum_{j=1}^{\infty} \frac{\alpha^{\lor} \otimes t_1^j}{j} \hspace{0.1cm}	u^j \Bigg)=\displaystyle\sum_{s=0}^{\infty}\hspace{0.1cm}p_\alpha^{(s)} \hspace{0.1cm}u^s,$$
where $p_\alpha^{(s)}$ denote the coefficient of $u^s$ in the power series $\textbf{P}_{\alpha}(u)$. Then clearly \hspace{0.02cm} $p_\alpha^{(0)}=1$. The following lemma is proved in \cite[Lemma 2.4]{khandai}.
\begin{lemma} \label{l1}
	For $\alpha \in R_{fin}^{+}$ and $ r\geq 1$ we have
	\begin{equation}\label{aq}
	(x_\alpha \otimes t_1)^r(x_{-\alpha} \otimes 1)^{r+1} - \displaystyle\sum_{s=0}^{r} (x_{-\alpha} \otimes t_1^{r-s})p_\alpha^{(s)} \in \hspace{0.2cm} \mathcal{U}(\mathcal{T}_{\mathrm{aff}}) (n_{\mathrm{aff}}^{+})
\end{equation} and
		\begin{equation}
	(x_\alpha \otimes t_1)^{r+1}(x_{-\alpha} \otimes 1)^{r+1} - p_\alpha^{(r+1)} \in \hspace{0.2cm} \mathcal{U}(\mathcal{T}_{\mathrm{aff}}) (n_{\mathrm{aff}}^{+}).
	\end{equation}
\end{lemma}

\vspace{0.3cm}
\subsection{Realisation of $\CT$ in terms of generators and relations}
\begin{definition}
	Let $\textbf{t}$ be the Lie algebra over $\C$ with the following presentation:
	
	{\bf Generators:}
	\begin{center}
	$\delta_{\underbar{\textbf{r}}}(\underbar{\textbf{s}})$,\hspace{0.2cm} $h_{i,\underbar{\textbf{k}}}$, \hspace{0.2cm} $e_{i,\underbar{\textbf{k}}}$, \hspace{0.2cm} $f_{i,\underbar{\textbf{k}}}$ , \hspace{0.2cm} $d_{j}$ \hspace{0.2cm} ($i=0,1,\cdots, l$,\hspace{0.1cm} $\underbar{\textbf{r}},\underbar{\textbf{s}},\underbar{\textbf{k}} \in \mathbb{Z}^{n-1},\hspace{0.2cm} 1\leq j \leq n$).\\
\end{center}

{\bf Relations:}
\begin{enumerate}
	\item[R1]
	\begin{enumerate}
		\item $\delta_{\underbar{\textbf{r}}}(\underbar{\textbf{s}})+\delta_{\underbar{\textbf{k}}}(\underbar{\textbf{s}})=\delta_{\underbar{\textbf{r}}+\underbar{\textbf{k}}}(\underbar{\textbf{s}})$;
		\vspace{0.1cm}
		\item
		$\delta_{\underbar{\textbf{r}}}(\underbar{\textbf{r}})=0$;
		\vspace{0.2cm}
		\item
		$[\delta_{\underbar{\textbf{r}}}(\underbar{\textbf{s}}), \delta_{\underbar{\textbf{p}}}(\underbar{\textbf{q}})]=[\delta_{\underbar{\textbf{r}}}(\underbar{\textbf{s}}), h_{i,\underbar{\textbf{k}}}]=[\delta_{\underbar{\textbf{r}}}(\underbar{\textbf{s}}), e_{i,\underbar{\textbf{k}}} ]=[\delta_{\underbar{\textbf{r}}}(\underbar{\textbf{s}}), f_{i,\underbar{\textbf{k}}} ]=0$;
		\vspace{0.2cm}
		\item 
		$[d_1, \delta_{\underbar{\textbf{r}}}(\underbar{\textbf{s}})]=0, \hspace{0.4cm} [d_j, \delta_{\underbar{\textbf{r}}}(\underbar{\textbf{s}})]=s_j \hspace{0.05cm} \delta_{\underbar{\textbf{r}}}(\underbar{\textbf{s}})$ for $j=2,3,\cdots, n$.
	\end{enumerate}
\vspace{0.2cm}
\item[R2] \hspace{0.2cm} $[h_{i,\underbar{\textbf{k}}} , \hspace{0.1cm}  h_{j,\underbar{\textbf{s}}}]=(h_i, h_j)\hspace{0.05cm} \delta_{\underbar{\textbf{k}}}(\underbar{\textbf{k}}+\underbar{\textbf{s}})$,
\vspace{0.2cm}
\item[R3]
\begin{enumerate}
	\item $[h_{i,\underbar{\textbf{k}}} ,\hspace{0.1cm}   e_{j,\underbar{\textbf{s}}}]=\alpha_j(h_i) \hspace{0.05cm} e_{j,\underbar{\textbf{k}}+\underbar{\textbf{s}}}$;
	\vspace{0.1cm}
	\item $[h_{i,\underbar{\textbf{k}}} ,\hspace{0.1cm}   f_{j,\underbar{\textbf{s}}}]=-\alpha_j(h_i)\hspace{0.05cm} f_{j,\underbar{\textbf{k}}+\underbar{\textbf{s}}}$,
\end{enumerate}
\vspace{0.2cm}
\item[R4] \hspace{0.2cm} $[e_{i,\underbar{\textbf{k}}} ,\hspace{0.1cm}   f_{j,\underbar{\textbf{s}}}]=\delta_{ij}\{h_{i, \underbar{\textbf{k}}+\underbar{\textbf{s}}}+\frac{2}{(\alpha_i, \alpha_i)}\delta_{\underbar{\textbf{k}}}(\underbar{\textbf{k}}+\underbar{\textbf{s}})\}$,
\vspace{0.2cm} 
\item[R5]
\begin{enumerate}
	\item $[e_{i,\underbar{\textbf{k}}} ,\hspace{0.1cm}   e_{i,\underbar{\textbf{s}}}]=0$;
	\item $[f_{i,\underbar{\textbf{k}}} ,\hspace{0.1cm}   f_{i,\underbar{\textbf{s}}}]=0$,
\end{enumerate}
\vspace{0.2cm}
\item[R6]
\begin{enumerate}
	\item (ad $e_{i,\underbar{\textbf{0}}})^{-A_{ij}+1}(e_{j,\underbar{\textbf{s}}})=0$,\hspace{1cm} $i\neq j$;
	\vspace{0.1cm}
	\item (ad $f_{i,\underbar{\textbf{0}}})^{-A_{ij}+1}(f_{j,\underbar{\textbf{s}}})=0$, \hspace{0.7cm} $i\neq j$ ,
\end{enumerate}
\vspace{0.2cm}
\item[R7]
\hspace{0.4cm} $[d_{j} ,\hspace{0.1cm}   e_{i,\underbar{\textbf{k}}}]=k_{j} \hspace{0.05cm}  e_{i,\underbar{\textbf{k}}}$,\hspace{0.2cm} $[d_{j} ,\hspace{0.1cm}   f_{i,\underbar{\textbf{k}}}]=k_{j} \hspace{0.05cm} f_{i,\underbar{\textbf{k}}}$,\hspace{0.2cm} $[d_{j} ,\hspace{0.1cm}   h_{i,\underbar{\textbf{k}}}]=k_{j} \hspace{0.05cm} h_{i,\underbar{\textbf{k}}}$, \hspace{0.4cm}($j=2,\cdots,n$)
	\vspace{0.2cm}
	\item[R8]
	\hspace{0.4cm} $[d_{1} ,\hspace{0.1cm}   e_{i,\underbar{\textbf{k}}}]=\delta_{i,0} \hspace{0.05cm}  e_{i,\underbar{\textbf{k}}}$,\hspace{0.2cm} $[d_{1} ,\hspace{0.1cm}   f_{i,\underbar{\textbf{k}}}]=-\delta_{i,0} \hspace{0.05cm} f_{i,\underbar{\textbf{k}}}$,\hspace{0.2cm} $[d_{1} ,\hspace{0.1cm}   h_{i,\underbar{\textbf{k}}}]=0$,
	\vspace{0.2cm}
	\item[R9]
	
	\hspace{0.5cm} $[d_i, \hspace{0.1cm} d_j]=0$ \hspace{1cm} ($1\leq i,j \leq n$).
\end{enumerate}
\end{definition}
\begin{theorem} [\cite{CMP}, Proposition 2.8]
	There is a Lie algebra isomorphism between $\textbf{t}$ and $\mathcal{T}$, the map is explicitly given by
	\vspace{0.2cm}
\begin{center}
	$e_{i,\underbar{\textbf{k}}} \mapsto e_i \otimes t^{\underbar{\textbf{k}}}$;\hspace{1cm}
	$f_{i,\underbar{\textbf{k}}} \mapsto f_i \otimes t^{\underbar{\textbf{k}}}$;\hspace{1cm}
	$h_{i,\underbar{\textbf{k}}} \mapsto h_i \otimes t^{\underbar{\textbf{k}}}$; \hspace{0.5cm}($i \in I$)\\
	\vspace{0.2cm}
	$e_{0,\underbar{\textbf{k}}} \mapsto f_\theta \otimes t_1t^{\underbar{\textbf{k}}}$;\hspace{0.5cm}
	$f_{0,\underbar{\textbf{k}}} \mapsto e_\theta \otimes t_1^{-1}t^{\underbar{\textbf{k}}}$; \hspace{0.4cm}
	$h_{0,\underbar{\textbf{k}}} \mapsto -h_{\theta} \otimes t^{\underbar{\textbf{k}}}+t^{\underbar{\textbf{k}}} K_1$;\\
	\vspace{0.3cm}
	$\delta_{\underbar{\textbf{r}}}(\underbar{\textbf{s}}) \mapsto \displaystyle\sum_{i=2}^n r_i t^{\underbar{\textbf{s}}} K_i$;\hspace{0.5cm}
	and \,\, $d_j \mapsto d_j$ $(j=1,2,\cdots, n).$
\end{center}
\end{theorem}
Let $\mathcal{T}^{+}$ be the Lie subalgebra of $\mathcal{T}$ generated by $\{e_{i,\underbar{\textbf{k}}}, f_{i,\underbar{\textbf{k}}}, d_1:  i \in I_{\mathrm{aff}} ,\hspace{0.1cm} \underbar{\textbf{k}} \in \Z^{n-1}_{\geq 0}\}$. Then we have
\begin{equation*}
	\begin{split}
		\mathcal{T}^{+} = &\mathfrak{g}\otimes \C[t_1^{\pm},t_2,\cdots, t_n] \hspace{0.2cm} \oplus \sum\limits_{\substack{i=1 \\ m_1 \in \Z,\hspace{0.1cm} m_i \geq 1\\ \underbar{\textbf{m}} \in \Z^{n-1}_{\geq 0}}}^n \C t_1^{m_{1}} t^{\underbar{\textbf{m}}} K_i \oplus \hspace{0.1cm} \C K_1 \oplus \hspace{0.1cm}\C d_1\\
		= &  \CT_{\mathrm{aff}}^{'} \otimes \C[t_2,t_3,\cdots, t_n] \hspace{0.2cm} \oplus \sum\limits_{\substack{i=2 \\ m_1 \in \Z,\hspace{0.1cm} m_i \geq 1\\ \underbar{\textbf{m}} \in \Z^{n-1}_{\geq 0}}}^n \C t_1^{m_{1}} t^{\underbar{\textbf{m}}} K_i \oplus \hspace{0.1cm} \hspace{0.1cm}\C d_1,
	\end{split}
	\end{equation*}
where $\CT_{\mathrm{aff}}^{'}$ is the derived algebra of $\CT_{\mathrm{aff}}$.
\begin{remark}
	Note that $t_1^{a_1}t_2^{a_2}\cdots t_j^{a_j}\cdots t_n^{a_n}K_j \in \mathcal{T}^{+}$ if and only if $ 0 \neq a_1 \in \mathbb{Z}, \,\, \underbar{\bf{a}} = (a_2, a_3, \cdots, a_n) \in \mathbb{Z}_{\geq 0}^{n-1}$ and ${\bf{a_j \geq 1}}$.
\end{remark}	
Let $\bar{\mathcal{T}}$ be the Lie subalgebra of $\mathcal{T}$ without $d_2, d_3,\cdots, d_n$:
$$\bar{\mathcal{T}}=\mathfrak{g} \otimes \C[t_1 ^{\pm},t_2^{\pm},\cdots, t_n^{\pm}] \hspace{0.2cm} \oplus \sum\limits_{\substack{i=1 \\ \textbf{m} \in \Z^{n}}}^n \C t^{\textbf{m}} K_i  \oplus \C d_1.$$ 

\subsection{Triangular decomposition:}
Let $N_{\mathcal{T}}^{+}$ and $N_{\mathcal{T}}^{-}$ be the Lie subalgebra of $\mathcal{T}$ generated by $\{e_{i,\underbar{\textbf{k}}}: \hspace{0.1cm} i=0,1,\cdots ,l, \hspace{0.2cm} \underbar{\textbf{k}} \in \Z^{n-1}\}$ and  $\{f_{i,\underbar{\textbf{k}}}: \hspace{0.1cm} i=0,1,\cdots ,l, \hspace{0.2cm} \underbar{\textbf{k}} \in \Z^{n-1}\}$ respectively. Set
\begin{center} 
	$N_{\mathcal{T}}^{0}= \mathfrak{h} \otimes \C[t_2^{\pm},t_3^{\pm},\cdots, t_n^{\pm}] \hspace{0.2cm} \oplus \sum\limits_{\substack{i=1 \\ \underbar{\textbf{m}} \in \Z^{n-1}}}^n \C t^{\underbar{\textbf{m}}} K_i$  $\oplus \hspace{0.1cm} \bigoplus\limits_{i=1}^n \C d_i.$ 
\end{center}
\begin{prop} We have the following:
	  	\begin{enumerate}
		\item   $N_{\mathcal{T}}^{+} = n_{\mathrm{aff}}^{+} \otimes \C[t_2^{\pm},t_3^{\pm},\cdots, t_n^{\pm}]  \hspace{0.2cm} \oplus \sum\limits_{\substack{i=1 \\ m_1 \geq 1,\\ \underbar{\textbf{m}} \in \Z^{n-1}}}^n \C t_1^{m_{1}} t^{\underbar{\textbf{m}}} K_i$
		\vspace{0.1cm}
		\item   $N_{\mathcal{T}}^{-} = n_{\mathrm{aff}}^{-} \otimes \C[t_2^{\pm},t_3^{\pm},\cdots, t_n^{\pm}]  \hspace{0.2cm} \oplus \sum\limits_{\substack{i=1 \\ m_1 \leq -1,\\ \underbar{\textbf{m}} \in \Z^{n-1}}}^n \C t_1^{m_{1}} t^{\underbar{\textbf{m}}} K_i$
		\vspace{0.2cm}
		\item $\mathcal{T} = N_{\mathcal{T}}^{-} \oplus N_{\mathcal{T}}^{0} \oplus N_{\mathcal{T}}^{+}$ is the triangular decomposition of $\mathcal{T}$.
	\end{enumerate}
\end{prop}
\vspace{0.3cm}
\par
  In $\mathcal{T}^{+}$, the elements $e_{i,\underbar{\textbf{k}}} \hspace{0.1cm} (i=0,1,\cdots ,l, \hspace{0.2cm} \underbar{\textbf{k}} \in \Z^{n-1}_{\geq 0})$ generate\\
  $$ N_{\mathcal{T}}^{+} \cap \mathcal{T}^{+} =  n_{\mathrm{aff}}^{+} \otimes \C[t_2,t_3,\cdots, t_n]  \hspace{0.2cm} \oplus \sum\limits_{\substack{i=1 \\ m_1 \geq 1, \hspace{0.1cm}m_i \geq 1 \\ \underbar{\textbf{m}} \in \Z^{n-1}_{\geq 0 }}}^n \C t_1^{m_{1}} t^{\underbar{\textbf{m}}} K_i,$$
  \par
  and the elements $f_{i,\underbar{\textbf{k}}} \hspace{0.1cm} (i=0,1,\cdots ,l, \hspace{0.2cm} \underbar{\textbf{k}} \in \Z^{n-1}_{\geq 0})$ generate\\
  $$ N_{\mathcal{T}}^{-} \cap \mathcal{T}^{+}= n_{\mathrm{aff}}^{-} \otimes \C[t_2,t_3,\cdots, t_n]  \hspace{0.2cm} \oplus \sum\limits_{\substack{i=1 \\ m_1 \leq -1, \hspace{0.1cm}m_i \geq 1\\ \underbar{\textbf{m}} \in \Z^{n-1}_{\geq 0}}}^n \C t_1^{m_{1}} t^{\underbar{\textbf{m}}} K_i.$$
  \par
  \subsection{Automorphisms of $\mathcal{T}$}\label{aut}
  Let $A = (a_{ij})$ be a $n \times n$ matrix in $GL_n{(\Z)}$ with determinant $\pm 1$. Then one has the following automorphism of $\mathcal{T}$ which we again denote by A:
  $$A(x \otimes t^{\textbf{m}}) = x \otimes t^{\textbf{m}A^{T}} \,\,\,\,\forall x \in \mathfrak{g}, \textbf{m} \in \Z^{n};$$
  $$A(t^{\textbf{m}}K_i) = \displaystyle{\sum_{r = 1}^{n}{a_{r i}t^{\textbf{m}A^{T}}K_r}};$$
  $$A(d_i) = \displaystyle{\sum_{r = 1}^n b_{ir}d_r},\,\, \mathrm{where}\,\, A^{-1} = (b_{ij}).$$
  
  \section{Global and local Weyl modules}\label{s2}
   In this section we define global and local Weyl modules for $\mathcal T$ . First we define two Lie algebra automorphisms $\Upsilon_0$ and $\Upsilon_\theta$ of $\mathcal{T}$ as given in \cite{ko};
   \vspace{0.1cm}
   $$\Upsilon_0=\mathrm{exp} \hspace{0.1cm} ad \hspace{0.1cm} e_0 \circ \mathrm{exp} \hspace{0.1cm} ad\hspace{0.1cm} (-f_0) \circ \mathrm{exp} \hspace{0.1cm} ad \hspace{0.1cm} e_0, $$ 
   $$\Upsilon_\theta=\mathrm{exp} \hspace{0.1cm} ad \hspace{0.1cm} e_\theta \circ \mathrm{exp}
    \hspace{0.1cm} ad\hspace{0.1cm} (-f_\theta) \circ \mathrm{exp} \hspace{0.1cm} ad \hspace{0.1cm} e_\theta. \vspace{0.2cm}$$ 
   Let M be a $\mathcal{T}$ module and if $M$ is integrable as a $\mathcal{T}_{\mathrm{aff}}$ module, then $\Upsilon_0,\Upsilon_{\theta} \in Aut \hspace{0.1cm}M$ are similarly defined. Moreover, they satisfy
   \vspace{0.2cm}
   \begin{center}
   $\Upsilon_0(xv)=\Upsilon_0(x) \Upsilon_0(v)$, \hspace{0.8cm} $\Upsilon_{\theta}(xv)=\Upsilon_{\theta}(x) \Upsilon_{\theta}(v)$ \hspace{0.8cm} for $x \in \mathcal{T},\hspace{0.1cm} v\in M.$
   \end{center}
\vspace{0.4cm}
\begin{lemma}[\cite{ko}, Lemma 2.8]\label{ab}
	$\Upsilon_0 \Upsilon_{\theta}\hspace{0.1cm}(e_{\theta} \otimes t_1^{m_{1}+2} t^{\underbar{\textbf{m}}})=\hspace{0.1cm} e_{\theta} \otimes t_1^{m_{1}} t^{\underbar{\textbf{m}}}$, \textrm{for} $m_1 \in \Z, \underbar{\textbf{m}} \in \Z^{n-1}.$ 
\end{lemma}
\subsection{Graded character}
Let $V$ be a $\mathcal{T}$-module with finite dimensional weight spaces. By restriction, $V$ is also a $\mathcal{T}^{+}$ and $\mathcal{T}_{\mathrm{aff}}$ module. As a ${\mathcal{T}_{\mathrm{aff}}}$-module we assume that $V$ has finite dimensional weight spaces
with respect to $\mathfrak{h}_{\mathrm{aff}}$. For any $\Lambda \in \mathfrak{h}^*$,  let $\Lambda - m \delta_1$ be a weight of the $\mathcal{T}_{\mathrm{aff}}$ module $V$. For a parameter $q_1$, the $q_1$-character of the $\mathcal{T}_{\mathrm{aff}}$-module $V$ is defined by $$\text{ch}_{q_1}(V) = \displaystyle{ \sum_{\Lambda \in \mathfrak{h}^*, m \in \mathbb{Z}}({\text{dim} \, V_{\Lambda -m \delta_1})\, e^{\Lambda}\, q_1^{m}}}.$$ 
Note that the Lie algebra $\mathcal{T}^{+}$ is $\mathbb{Z}^{n-1}_{\geq 0}$-graded with the grading:
$\text{deg}(a \otimes t_1^{r_1}t_2^{r_2}\cdots t_n^{r_n}) : = (r_2, \ldots, r_n) $ for $a \in \mathfrak{g}$, $\text{deg} (t_1^{s_1}t_2^{s_2}\cdots t_n^{s_n} K_i) := (s_2, \ldots, s_n)$ for $1 \leq i \leq n$ and $\text{deg} \,d_1 = 0$. Each weight space
$V_{\Lambda -m \delta_1}$ has graded decomposition $V_{\Lambda -m \delta_1} := \displaystyle{\sum_{\underbar{\textbf{r}} \in \mathbb{Z}^{n-1}}}{V_{\Lambda -m \delta_1} [\underbar{\textbf{r}}]}$. For parameters $q_1,q_2, \cdots ,q_n$, the graded character of the module $V$ is defined as
$$\text{ch}(V)_{q_1, \ldots, q_n} := \sum_{\substack{\Lambda \in \mathfrak{h}^*, m \in \mathbb{Z}\\\underbar{\textbf{r}} \in \mathbb{Z}^{n-1} }}{\text{dim}\,V_{\Lambda -m \delta_1}[\underbar{\textbf{r}}] e^{\Lambda} \, q_1^{m} q_2^{r_2}\cdots q_{n}^{r_n}}.$$

\begin{definition}[\textbf{Global Weyl module}]
	 Let $\Lambda$ be a dominant integral weight of $\mathcal{T}_{\mathrm{aff}}$. Then the global Weyl module $W_{\mathrm{glob}}{(\Lambda)}$ is a highest weight module for $\mathcal{T}$ with highest weight $\Lambda$ and is generated by $v_\Lambda$ with the following defining relations:
	 \vspace{0.3cm}
	  \begin{center} 
	  	$e_{i,\underbar{\textbf{k}}}.\hspace{0.05cm} v_{\Lambda}=0$ \hspace{0.2cm}($i \in I_{\mathrm{aff}}, \hspace{0.1cm} \underbar{\textbf{k}} \in \Z^{n-1}$),\hspace{0.8cm} $h.v_{\Lambda}= \Lambda(h)\hspace{0.05cm}v_{\Lambda}$ \hspace{0.4cm}($h \in \mathfrak{h}_{\mathrm{aff}}$)\\
	  	\vspace{0.3cm}
	  	$\hspace{0.8cm}f_{i}^{\Lambda(h_i)+1}. \hspace{0.05cm}v_{\Lambda}=0$ \hspace{0.4cm} ($i \in I_{\mathrm{aff}}$), $\hspace{0.8cm} d_{i}.\hspace{0.05cm}v_{\Lambda}=0$ \hspace{0.2cm}($i=2,3,\cdots,n$)\\
		\vspace{0.3cm}
		and $t^{\underbar{\textbf{m}}} K_i  v_{\Lambda} = 0$ for $\underbar{\textbf{m}} \in \mathbb{Z}^{n-1}$  \hspace{0.2cm}($i=2,3,\cdots,n$).

	  	\end{center} 
  	\par
  	\vspace{0.2cm}
  	The global Weyl module $W^+_{\mathrm{glob}}{(\Lambda)}$ is a highest weight module for $\mathcal{T}^{+}$ with highest weight $\Lambda$ and is generated by $v_\Lambda^+$ with the following defining relations:
  	\vspace{0.2cm}
  	\begin{center} 
  		$e_{i,\underbar{\textbf{k}}}.\hspace{0.05cm} v_{\Lambda}^{+}=0$ \hspace{0.2cm}($i \in I_{\mathrm{aff}}, \hspace{0.1cm} \underbar{\textbf{k}} \in \Z^{n-1}_{\geq 0}$),\hspace{0.8cm} $h.v_{\Lambda}^{+}= \Lambda(h)\hspace{0.05cm}v_{\Lambda}^{+}$ \hspace{0.4cm}($h \in \mathfrak{h}_{\mathrm{aff}}$)\\
  		\vspace{0.25cm}
  		$f_{i}^{\Lambda(h_i)+1}. \hspace{0.05cm}v_{\Lambda}^{+}=0$ \hspace{0.7cm} ($i \in I_{\mathrm{aff}}$)\\
  		\vspace{0.25cm} and $t^{\underbar{\textbf{m}}} K_i  v_{\Lambda} = 0$ for $\underbar{\textbf{m}} \in \mathbb{Z}^{n-1}_{\geq 0}, m_i \geq 1$  \hspace{0.2cm}($i=2,3,\cdots,n$).\\
  	\end{center}
\end{definition}

Let 
$$\bar{ \mathfrak{h}}_{\mathcal{T}} := N_{\mathcal{T}}^0 \cap \bar{\mathcal{T}} = \mathfrak{h} \otimes \mathbb{C}[t_{2}^{\pm 1}, \cdots, t_n^{\pm 1}]  \hspace{0.2cm} \oplus \sum\limits_{\substack{i=1 \\ \underbar{\textbf{m}} \in \Z^{n-1}}}^n \C t^{\underbar{\textbf{m}}} K_i \oplus \mathbb{C} d_1$$
$$=\mathfrak{h}_{\mathrm{aff}}^{'} \otimes \mathbb{C}[t_{2}^{\pm}, \cdots, t_n^{\pm}]  \hspace{0.2cm} \oplus \sum\limits_{\substack{i=2\\ \underbar{\textbf{m}} \in \Z^{n-1}}\,}^n \C t^{\underbar{\textbf{m}}} K_i \oplus \mathbb{C} d_1,$$
where $\mathfrak{h}_{\mathrm{aff}}^{'} := \mathfrak{h} \oplus \mathbb{C}K_1$.  In the last expression we identify $t^{\underbar{\textbf{m}}} K_1$ with $K_1 \otimes t^{\underbar{\textbf{m}}} \,\, \forall \,\, \underbar{\textbf{m}} \in \mathbb{Z}^{n-1}$.

 \subsection{The algebra $A(\Lambda)$}
 Following the approach of \cite{CK}, $W_{\mathrm{glob}}{(\Lambda)}$ can be given a right $ \bar{ \mathfrak{h}}_{\mathcal{T}}$-module structure: $X.(v_{\Lambda}).a = X(a.v_{\Lambda})$, where $X \in Y \in \mathcal{U}(\bar{\mathcal{T}})$ and $a \in \bar{ \mathfrak{h}}_{\mathcal{T}}$. Argument similar
 to (\cite{CK}, 3.4) implies that this right action is well defined.
 We note that the weight spaces of $W_{\mathrm{glob}}(\Lambda)$ are invariant under the right $ \bar{ \mathfrak{h}}_{\mathcal{T}}$ action. In particular $W_{\mathrm{glob}}(\Lambda)_{\Lambda}$ is invariant under $ \bar{ \mathfrak{h}}_{\mathcal{T}}$. Now consider the space $A(\Lambda) := \frac{\mathcal{U}( \bar{ \mathfrak{h}}_{\mathcal{T}})}{\mathrm{Ker}_{\mathcal{U} (\bar{ \mathfrak{h}}_{\mathcal{T}})} v_{\Lambda}}$, where $\mathrm{Ker}_{\mathcal{U} (\bar{ \mathfrak{h}}_{\mathcal{T}})} v_{\Lambda} = \{ x \in \mathcal{U} (\bar{ \mathfrak{h}}_{\mathcal{T}}): v_{\Lambda}.x = x.v_{\Lambda} = 0\}$. Then the 
 assignment $a \mapsto v_{\Lambda}.a = a.v_{\Lambda}$ gives a vector space isomorphism between  $A(\Lambda)$ and $W_{\mathrm{glob}}(\Lambda)_{\Lambda}$ . 
 
 Let $\Lambda = \sum_{i = 0}^{l}{r_i \Lambda_i}$. Observe that in this case $A(\Lambda)$ is a commutative subalgebra of $\bar{ \mathfrak{h}}_{\mathcal{T}}$ as $t^{\underbar{\textbf{m}}} K_i v_{\Lambda} = 0 = d_1.v_{\Lambda} \,\, \forall \,\, \underbar{\textbf{m}} \in \mathbb{Z}^{n-1}, 2 \leq i \leq n$. Let $L:= \mathbb{C}[t_2^{\pm}, \ldots, t_n^{\pm}]$. Consider the algebra $L^{\otimes N}$, where $N = \sum_{i = 0}^{l}{r_i}$.  Let $\mathcal{S}_{n}$ denote the symmetric group on $n$-letters. Define $\mathcal{S}_{r_{\Lambda}} = \mathcal{S}_{r_0} \times \cdots \times \mathcal{S}_{r_l}$ and $B(\Lambda) := (L^{\otimes N})^{\mathcal{S}_{r_\Lambda}} :=( L^{\otimes r_0})^{\mathcal{S}_{r_0}} \otimes \cdots  \otimes ( L^{\otimes r_l})^{\mathcal{S}_{r_l}}$. Then we have the following:

\begin{thm} \label{t1}
$A(\Lambda) \cong B(\Lambda)$ as algebras.
\end{thm}
To prove the theorem we use identities by Garland in the universal enveloping algebra of $\mathcal{T}_{\mathrm{aff}}$ stated in Lemma \ref{l1}. 
Now for $s \in \mathbb{Z}$, consider the $\mathfrak{sl}_2$ triple $\{x_{\alpha} \otimes t_1^{s}, x_{-\alpha} \otimes t_1^{-s}, \alpha^{\vee}- \frac{2s}{(\alpha,\alpha)}K_1\}$. For a monomial $a \in L$, we have a Lie algebra isomorphism
$\phi_{\alpha, s}^{a}: \mathfrak{sl}_2 \otimes \mathbb{C}[t_1^{\pm 1}] \rightarrow \mathfrak{sl}_2 \otimes \mathbb{C}[a^{\pm 1}]$ defined by $\phi_{\alpha, s}^{a}(x_{\pm \alpha} \otimes t_1^{p}) = x_{\pm \alpha} \otimes t_1^{\pm s} \otimes a^p$ and
$\phi_{\alpha, s}^{a}(\alpha^{\vee} \otimes t_1^{p}) = (\alpha^{\vee}- \frac{2s}{(\alpha,\alpha)}K_1) \otimes a^{p}$. Using this isomorphism, for $\gamma = \alpha + s \delta_1 \in R_{\mathrm{aff}}^{re^+}$, we have the following Garland identities:
$$(x_\gamma \otimes a)^r(x_{-\gamma} \otimes 1)^{r+1} v_{\Lambda} = \displaystyle\sum_{m=0}^{r} (x_{-\gamma} \otimes a^{r-m})p_{\gamma, a}^{(m)} v_{\Lambda}$$
	and  $$(x_\gamma \otimes a)^{r+1}(x_{-\gamma} \otimes 1)^{r+1} v_{\Lambda} = p_{\gamma, a}^{(r+1)} v_{\Lambda},$$ where $p_{\gamma, a }^{(s)} := \phi_{\alpha, s}^{a}(p_{\alpha}^{(s)})$ and $a \in L$ is any monomial.
	  We also have $$x_{\gamma}(x_\gamma \otimes a)^r(x_{-\gamma} \otimes 1)^{r+1} v_{\Lambda} = \sum_{m = 0}^{r}{ (\gamma^{\vee} \otimes a^{r - m} )p_{\gamma, a}^{(m)}} \,v_{\Lambda}.$$
	  
Now, as an application of PBW theorem we have that $I_{\Lambda} := \mathrm{Ker}_{\mathcal{U} (\bar{ \mathfrak{h}}_{\mathcal{T}})}.v_{\Lambda}$ is an ideal of $\mathcal{U} (\bar{ \mathfrak{h}}_{\mathcal{T}})$ generated by the following elements:
$$h - \Lambda(h) \hspace{0.5cm} (h \in \mathfrak{h}_{\mathrm{aff}}), \hspace{1cm}  t^{\underbar{\textbf{m}}} K_j \hspace{0.5cm} (\, {\underbar{\textbf{m}}} \in \mathbb{Z}^{n-1}, 2 \leq j \leq n),\hspace{0.5cm} d_1, \hspace{0.5cm} p_{\alpha_i, a}^{(r)}$$
$$\mathrm{and} \,\,\, \sum_{m = 0}^{r-1}{ (\alpha_i^{\vee} \otimes a^{r - m} )p_{\gamma, a}^{(m)}},\,\,\,\,\,\,\, \,\,\,\,
 \forall \,\, \mathrm{monomials} \,\, a \,\, \mathrm{in} \,\, L, \, |r| > \Lambda(\alpha_i^{\vee}), \,\mathrm{and}\,\,  i \in I_{\mathrm{aff}}.$$

 Define a map
 $\Phi : \mathcal{U} (\bar{ \mathfrak{h}}_{\mathcal{T}}) \rightarrow B(\Lambda)$ by $$\Phi(\alpha_{i}^{\vee} \otimes a) = 1^{\otimes (r_0+ \cdots r_{i-1})} \otimes \big(\sum_{k = 0}^{r_i - 1}{1^{\otimes k} \otimes a \otimes 1^{\otimes r_i - k -1}}\big) \otimes 1^{\otimes(r_{i+1} + \cdots  + r_n)}, \,\,\,0 \leq i \leq l, a \in L$$ and $\Phi(d_1) = 0$.
 Then note that $\Phi$ is an algebra homomorphism which maps the elements $t^{\bf{m}}K_i$ to zero for $2 \leq i \leq n$. We then have the following proposition:
 \begin{prop} \label{imp}
 The map $\Phi$ factors through $I_{\Lambda}$ and the induced map $\tilde{\Phi}: A ({\Lambda}) \rightarrow B(\Lambda)$ is an algebra isomorphism.
 \end{prop}
 \begin{proof}
 Clearly $\Phi$ is a surjective algebra homomorphism. It is easy to see that $\Phi(h - \Lambda(h)) = 0$ for all $h \in \mathfrak{h}_{\mathrm{aff}}$. We observe that $\Phi$ maps the element $\Phi(\alpha_{i}^{\vee} \otimes a)$ to a power sum symmetric function. The rest of the generators of $I_{\Lambda}$ 
 are mapped to zero by definition of $\Phi$  and identities between power sum symmetric functions and elementary symmetric functions (see Appendix at the end of this paper for the identities \ref{e1} and \ref{e2}). So we have  a surjective algebra homomorphism $\tilde{\Phi}: A ({\Lambda}) \rightarrow B(\Lambda)$.  Argument for injectivity of $\tilde{\Phi}$ is similar to that of proof of Lemma 10 of \cite{CK}. 
 \end{proof}
 
 From the above result it follows that $A(\Lambda)$ is finitely generated as an algebra. Now, repeating the same argument for $\CT^{+}$, we have that $A(\Lambda)^{+}: = \frac{\mathcal{U}( \bar{ \mathfrak{h}}_{\mathcal{T}^+})}{\mathrm{Ker}\,_{\mathcal{U}(\bar{\fh}_{\CT^{+}})} v_{\Lambda}^{+}}$ is isomorphic to $B(\Lambda)^{+}$ as an algebra, where  
 \begin{equation*}
 	\begin{split}
\bar{\fh}_{\CT^{+}} &= \mathfrak{h} \otimes \mathbb{C}[t_{2}, \cdots, t_n]  \hspace{0.2cm} \oplus \sum\limits_{\substack{\hspace{0.8cm}i=1  \\ \underbar{\textbf{m}} \in \Z^{n-1}_{\geq 0}},\, m_i \geq 1}^n \C t^{\underbar{\textbf{m}}} K_i \oplus \mathbb{C} d_1\\
&= \mathfrak{h}_{\mathrm{aff}}^{'} \otimes \mathbb{C}[t_{2}, \cdots, t_n]  \hspace{0.2cm} \oplus \sum\limits_{\substack{\hspace{0.8cm}i=2  \\ \underbar{\textbf{m}} \in \Z^{n-1}_{\geq 0}},\, m_i \geq 1}^n \C t^{\underbar{\textbf{m}}} K_i \oplus \mathbb{C} d_1,
	\end{split}
\end{equation*}
and 
$$A(\Lambda)^+ = ({(L^+)}^{\otimes N})^{\mathcal{S}_{r_\Lambda}} =( {(L^+)}^{\otimes r_0})^{\mathcal{S}_{r_0}} \otimes \cdots  \otimes ( {(L^+)}^{\otimes r_l})^{\mathcal{S}_{r_l}}$$ with ${L^+} := \mathbb{C}[t_{2}, \cdots, t_n] $.

\subsection{Local Weyl modules for $\CT$ and $\CT^+$}
For a dominant weight $\Lambda$ of $\CT_{\mathrm{aff}}$, let $M$ be a maximal ideal of $A(\Lambda)$ and $\mathbb{C}_{\Lambda, M}$ be the corresponding one dimensional representation of $A(\Lambda)$. Then $W_{\mathrm{loc}}{(\Lambda, M)}:= W_{\mathrm{glob}}(\Lambda) \otimes_{A(\Lambda)}\mathbb{C}_{\Lambda, M} $
is called a local Weyl module for $\mathcal{T}$. Similarly define ${W^+}_{\mathrm{loc}}{(\Lambda, M^+)}:= {W^+}_{\mathrm{glob}}(\Lambda) \otimes_{A(\Lambda)^+}\mathbb{C}_{\Lambda, M^+},$ where $M^+$ is a maximal ideal of $L^+$.

\begin{lemma}\label{red}  Let $\Lambda$ be a dominant integral weight of $\mathcal{T}_{\mathrm{aff}}.$
 
		$(\mathrm{i})$ For any positive root $\gamma \in R_{\mathrm{aff}}^{+}$, there exists $N(\gamma) \in  \Z_{+}$ such that for any $\underbar{\textbf{a}}=(a_2,\cdots,a_n) \in \Z^{n-1}$ we have
		\vspace{0.1cm} 
		$$(x_{-\gamma} \otimes \hspace{0.05cm} t^{\underbar{\textbf{a}}}) \hspace{0.05cm} .v_{\Lambda} \hspace{0.1cm} \in \displaystyle\sum_{0\leq m_i \leq N(\gamma)} (x_{-\gamma} \otimes \hspace{0.05cm}t_2^{m_2}t_3^{m_3}\cdots t_n^{m_n})\, A(\Lambda) \hspace{0.05cm}v_{\Lambda}.$$
		
	$(\mathrm{ii})$ For any positive integer $q>0$, there exists $N(q) \in  \Z_{\geq 0}$ such that for any $\underbar{\textbf{a}}=(a_2,\cdots,a_n) \in \Z^{n-1}$ and $i=2,\cdots,n$ we have 
		\vspace{0.1cm}
		$$t_1^{-q}t^{\underbar{\textbf{a}}} K_{i} \hspace{0.05cm} .v_{\Lambda} \in \hspace{-0.25cm} \displaystyle\sum_{0\leq m_j \leq N(q)} \hspace{-0.35cm} t_1^{-q}t_2^{m_2}\cdots t_n^{m_n} K_i \hspace{0.1cm} A(\Lambda) \hspace{0.05cm}v_{\Lambda} + \hspace{-0.25cm} \displaystyle\sum_{0\leq m_j \leq N(q)} \hspace{-0.35cm} {(\mathcal{T}_{\mathrm{aff}})}_{-q\delta_1} \hspace{0.1cm}\otimes \hspace{0.05cm}t_2^{m_2}\cdots t_n^{m_n} A(\Lambda) \hspace{0.05cm}v_{\Lambda}.$$
\end{lemma}
\begin{proof}

	Let $\gamma=\alpha+r\delta_1$ where $r>0, \alpha \neq 0$, $a=t_3^{a_3}t_4^{a_4}\cdots t_n^{a_n} \in \C[t_1^{\pm},t_2^{\pm},\cdots, t_n^{\pm}]$ and let  $\gamma_{0}:=[x_\alpha \otimes t_1^{r} a^{-1},\hspace{0.1cm} x_{-\alpha} \otimes t_1^{-r} a]$. Let $\mathcal{S}$ be a subalgebra of \,$\mathcal{T}$ generated by $\{x_\alpha \otimes t_1^{r} a^{-1},\hspace{0.1cm} x_{-\alpha} \otimes t_1^{-r} a,\gamma_0 \}$. Then it is easy to see that $\mathcal{S} $ is isomorphic to $\mathfrak{sl}_2(\C) \cong <\{x_\alpha, x_{-\alpha}, \alpha^{\lor}\}>$.
	Let $\Phi_2: \mathfrak{sl}_2(\C) \otimes \C[t_2^{\pm 1}] \to \mathcal{S} \otimes \C[t_2^{\pm 1}]$ be defined by
	$$x_\alpha \otimes t_2^{p} \mapsto x_\alpha \otimes t_1^{r} a^{-1} \otimes t_2^{p} ,  \,\, \,\, x_{-\alpha} \otimes t_2^{p} \mapsto x_{-\alpha} \otimes t_1^{-r} a \otimes t_2^{p}  \,\, and \,\, \alpha^{\lor} \otimes t_2^{p} \mapsto \gamma_{0} \otimes t_2^{p}.$$ Then $\Phi_2$ is an isomorphism. 
	Let $\Phi_2(p_{\alpha}^{(m_2)})=p_{\alpha,\Phi_2 }^{(m_2)}$. Then using \eqref{aq}, we get  that $$\displaystyle\sum_{m_2=0}^{\gamma(\alpha^\lor)} x_{-\alpha} \otimes t_1^{-r}at_2^{\gamma(\alpha^\lor)-m_2}\hspace{0.1cm} \Phi_2(p_{\alpha}^{(m_2)}) v_{\Lambda}=0.$$ 
	
	From this we get that 
$$x_{-\gamma} \otimes t_2^{\gamma(\alpha^\lor)}t_3^{a_3}\cdots t_n^{a_n}v_{\Lambda}=-\displaystyle\sum_{m_2=1}^{\gamma(\alpha^\lor)} x_{-\gamma} \otimes t_2^{\gamma(\alpha^\lor)-m_2}t_3^{a_3}\cdots t_n^{a_n}\hspace{0.1cm} p_{\alpha,\Phi_2 }^{(m_2)}v_{\Lambda}.$$
	Now for $\gamma=r\delta_{1} (r>0)$, the root vectors corresponding to $-\gamma$ in $\mathcal{T}_{\mathrm{aff}}$ are $\{h_i \otimes t_1^{-r}: 1 \leq i \leq l \}$. Proceeding as above for $\gamma= \alpha_{i}+ r\delta_1$ and applying $x_{\alpha_i}$, we get that
	$$\displaystyle\sum_{m_2=0}^{\gamma(h_i)} h_{i} \otimes t_1^{-r}t_2^{\gamma(h_i)-m_2}t_3^{a_3}\cdots t_n^{a_n}\hspace{0.1cm} p_{\alpha,\Phi_2 }^{(m_2)}v_{\Lambda}=0.$$
	Using the techniques developed in \cite[Corollary 3.1]{2002}, we get that for any positive root $\gamma \in R_{\mathrm{aff}}^{+}$, $x_{-\gamma} \otimes \hspace{0.05cm} t_2^{a_2} t_3^{a_3} \cdots t_n^{a_n}$ lies in the span of $\{x_{-\gamma} \otimes \hspace{0.05cm} t_2^{m_2} t_3^{a_3} \cdots t_n^{a_n}p_{\alpha, \Phi_2}^{(N(\gamma)-m_2)} v_{\Lambda}: 0 \leq m_2 \leq N(\gamma)\}$, where
	\[
	N(\gamma) =
	\begin{cases}
		\gamma(\alpha^\lor)& \text{if } \gamma= \alpha +r \delta_1, \alpha \neq 0,\\
		max \,                                                                                                                                                                                                                                                                                                                                                                                                                                                                                                                                                                                                                                                                                                                                                                                                                                                                                                                                                                                                                                                                                                                                                                                                                                                                                             \{ \gamma(h_i):1 \leq i \leq l\}              & \text{if } \gamma= r \delta_1.\\
	\end{cases}
\] Now
	\begin{equation}\label{2n}
		\begin{aligned}
				&x_{-\gamma} \otimes \hspace{0.05cm} t_2^{m_2} t_3^{a_3} \cdots t_n^{a_n}p_{\alpha, \Phi_2}^{(N(\gamma)-m_2)} v_{\Lambda}\\ 
				&\hspace{-0.45cm}=p_{\alpha, \Phi_2}^{(N(\gamma)-m_2)} \hspace{0.1cm} (x_{-\gamma} \otimes \hspace{0.05cm} t_2^{m_2} t_3^{a_3} \cdots t_n^{a_n}) \hspace{0.05cm} v_{\Lambda} + \hspace{0.05cm} [ \hspace{0.05cm} x_{-\gamma} \otimes \hspace{0.05cm} t_2^{m_2} t_3^{a_3} \cdots t_n^{a_n}, \hspace{0.05cm} p_{\alpha, \Phi_2}^{(N(\gamma)-m_2)}] \hspace{0.1cm} v_{\Lambda}.
		\end{aligned}
	\end{equation}
Defining an appropriate isomorphism $\Phi_3$ similar to $\Phi_2$ as defined above, the 1st term of \eqref{2n} is equal to
\begin{equation}\label{3}
		\begin{aligned}
&p_{\alpha, \Phi_2}^{(N(\gamma)-m_2)} \displaystyle\sum_{m_3=0}^{N(\gamma)} \lambda_{m_3} \hspace{0.1cm} x_{-\gamma} \otimes \hspace{0.05cm} t_2^{m_2} t_3^{m_3} \cdots t_n^{a_n}\hspace{0.05cm} p_{\alpha, \Phi_3}^{(N(\gamma)-m_3)} v_{\Lambda}\\
&\hspace{-0.5cm}=\displaystyle\sum_{m_3=0}^{N(\gamma)} \lambda_{m_3} \hspace{0.1cm} x_{-\gamma} \otimes \hspace{0.05cm} t_2^{m_2} t_3^{m_3} t_4^{a_4} \cdots t_n^{a_n} \, p_{\alpha, \Phi_2}^{(N(\gamma)-m_2)} p_{\alpha, \Phi_3}^{(N(\gamma)-m_3)} \hspace{0.05cm} v_{\Lambda}\\    &\hspace{1.5cm}+ \displaystyle\sum_{m_3=0}^{N(\gamma)} \lambda_{m_3} \hspace{0.05cm} [\hspace{0.05cm}p_{\alpha, \Phi_2}^{(N(\gamma)-m_2)},\hspace{0.05cm}  x_{-\gamma} \otimes \hspace{0.05cm} t_2^{m_2} t_3^{m_3} t_4^{a_4}  \cdots t_n^{a_n}] \hspace{0.05cm} p_{\alpha, \Phi_3}^{(N(\gamma)-m_3)} v_{\Lambda}.
\end{aligned}
 \end{equation}

	Note that the coefficient of $t_2$ in $p_{\alpha, \Phi_{2}}^{(s)}$ is always less than or equal to $s$. Hence the second term of \eqref{3} is in the span of $\{x_{-\gamma} \otimes \hspace{0.05cm} t_2^{m_2} t_3^{m_3} \cdots t_n^{a_n}p_{\alpha, \Phi_3}^{(N(\gamma)-m_3)} v_{\Lambda}: 0 \leq m_2, m_3 \leq N(\gamma)\}$ and second term of \eqref{2n} belongs to the span of $\{x_{-\gamma} \otimes \hspace{0.05cm} t_2^{m_2} t_3^{a_3} \cdots t_n^{a_n} v_{\Lambda}: 0 \leq m_2 \leq N(\gamma)\}$. Therefore, $x_{-\gamma} \otimes \hspace{0.05cm} t_2^{a_2} t_3^{a_3} \cdots t_n^{a_n} v_{\Lambda}$ lies in the span of $\{x_{-\gamma} \otimes \hspace{0.05cm} t_2^{m_2} t_3^{m_3} t_4^{a_4} \cdots t_n^{a_n} A(\Lambda) v_{\Lambda}: 0 \leq m_2, m_3 \leq N(\gamma)\}$.
	Using similar argument as above, we get that 
	$$(x_{-\gamma} \otimes \hspace{0.05cm} t_2^{a_2} t_3^{a_3} \cdots t_n^{a_n}) \hspace{0.05cm} .v_{\Lambda} \hspace{0.1cm} \in \displaystyle\sum_{0\leq m_i \leq N(\gamma)} (x_{-\gamma} \otimes \hspace{0.05cm}t_2^{m_2}t_3^{m_3}\cdots t_n^{m_n}) A(\Lambda) \hspace{0.05cm}v_{\Lambda}.$$
	That completes the proof of the assertion (i). For (ii) we note that 
	$$(t_1^{-q}t_2^{a_2} \cdots t_i^{a_i}\cdots t_n^{a_n} K_{i}) \hspace{0.05cm} v_{\Lambda}= ([x_{\alpha} \otimes t_i, x_{-\alpha} \otimes t_1^{-q}t_2^{a_2} \cdots t_i^ {a_i-1} \cdots t_n^{a_n}] -h_{\alpha} \otimes t_1^{-q}t_2^{a_2} \cdots t_n^{a_n} \hspace{0.05cm}) v_{\Lambda}.$$ Then the proof follows by using the technique similar to \cite[Lemma 3.10(ii)]{ko}.
\end{proof}

\begin{cor} \label{c1}
Let $\Lambda$ be a dominant integral weight of $\mathcal{T}_{\mathrm{aff}}$ and let $\gamma \in R^{+}_{\textrm{aff}}$. Then the weight space ${W_{\mathrm{glob}}(\Lambda)}_{\Lambda - \gamma}$ is finitely generated over $A(\Lambda)$. In particular, we have ${W_{\mathrm{loc}}{(\Lambda, M)}}_{\Lambda - \gamma}$ is finite dimensional and
$W_{\mathrm{loc}}{(\Lambda, M)} = \mathcal{U}(\CT^{+}) v_{\Lambda, M}$, where $W_{\mathrm{loc}}{(\Lambda, M)}=W_{\mathrm{glob}}(\Lambda)\otimes_{A(\Lambda)} \mathbb{C}v_{\Lambda, M}$. 
\end{cor}

\begin{proof}
The proof follows from Lemma \ref{red}.
\end{proof} 
\section{Level one case} \label{s3}	
In this section we use the results in Section \ref{s2} to the special case where $\Lambda = \Lambda_0$ and obtain an upper bound for the graded character of $W_{\mathrm{loc}}{(\Lambda_0, M)}$. 
In this case we have $\tilde{\Phi}(A(\Lambda_0)) = L_0$, a Laurent polynomial ring in $n-1$ variables and we denote it by 
$L_0: = \mathbb{C}[x_2^{\pm1}, \ldots, x_{n}^{\pm1}]$ and under this isomorphism we identify $A(\Lambda_0)$ with $L_0$. Then $W_{\mathrm{glob}}(\Lambda_0)$ is a right $L_0$ module ;
$ (Y .v_{\Lambda_0}).x^{\underbar{\textbf{k}}} = Y. ({\tilde{\Phi}}^{-1}(x^{\underbar{\textbf{k}}}).v_{\Lambda_0})) = Y.(h_{0, \underbar{\textbf{k}}}.v_{\Lambda_0})$, where $Y \in \mathcal{U}(\bar{\mathcal{T}})$ and $(k_2, \cdots, k_n) =\underbar{\textbf{k}} \in \mathbb{Z}^{n-1}$
, $x^{\underbar{\textbf{k}}} = x_2^{k_2} \cdots x_n^{k_n} \in L_0$. We have the following observations:
\begin{itemize}
\item $\forall\,\,\, \underbar{\textbf{k}} \in \mathbb{Z}^{n-1}$, \,\,$h_{i, {\underbar{\textbf{k}}}} v_{\Lambda_0} = 0$ for all $i \in \{1,2, \ldots l\}$ and
\item $\forall\,\,\, \underbar{\textbf{k}} \in \mathbb{Z}^{n-1}$, $f_{i, \underbar{\textbf{k}}} v_{\Lambda_0} = 0$ for all $i \in \{1,2, \ldots l\}$.
\end{itemize}

We have $W_{\mathrm{loc}}{(\Lambda_0, M)} = W_{\mathrm{glob}}(\Lambda_0) \otimes_{A(\Lambda_0)} \mathbb{C}_M$, where $M$ is a maximal ideal of $L_0$ and $\mathbb{C}_M$ is a one dimensional representation of $A(\Lambda_0) = L_0$.
By Corollary \ref{c1}, it follows that $W_{\mathrm{loc}}{(\Lambda_0, M)} = \mathcal{U}(\mathcal{T}^+) \otimes_{A(\Lambda_0)} {\mathbb{C}v_{\Lambda_0, M}}$. Similarly, for the subalgebra $\CT^+$ we have $W^{+}_{\mathrm{loc}}{(\Lambda, M^+}) = \mathcal{U}(\mathcal{T}^+) \otimes_{{A(\Lambda_0)}^+} {\mathbb{C}{v_{\Lambda_0, M^+}^+}}$ where ${A(\Lambda_0)}^+$ is identified with ${L^+_0} := \mathbb{C}[x_2, \ldots, x_{n}]$.  As maximal ideals of $L_0^{+}$ are
in a one to one correspondence with the points in $\mathbb{C}^{n-1}$, we  write $W^+_{\mathrm{loc}}{(\Lambda_0, M^+)} = W^+_{\mathrm{loc}}{(\Lambda_0, \underbar{\bf{a}})}$ for some $\underbar{\bf{a}} \in \mathbb{C}^{n-1}$. In this case we have $v^+_{\Lambda_0, \underbar{\bf{a}}}.x^{\underbar{\textbf{k}}} = h_{0, \underbar{\textbf{k}}}\, v^+_{\Lambda_0, \underbar{\bf{a}}} = \underbar{\bf{a}}^{\underbar{\textbf{k}}}v^+_{\Lambda_0, \underbar{\bf{a}}} \,\, \forall \,\,\underbar{\textbf{k}} \in \mathbb{Z}^{n-1}_{\geq 0}$, where $v^+_{\Lambda_0, \underbar{\bf{a}}} := v^+_{\Lambda_0, M^+}$. Following similar argument as in Proposition 3.14 of \cite{ko} we have
$$\mathrm{ch}_{q_1} \big(W^+_{\mathrm{loc}}{(\Lambda_0, \underbar{\bf{a}})}\big) = \mathrm{ch}_{q_1}\big(W^+_{\mathrm{loc}}{(\Lambda_0, \underbar{\bf{0}})}\big)\,\, \forall \,\,\underbar{\bf{a}} \in \mathbb{C}^{n-1}.$$
In the rest of this section we work with $W^+_{\mathrm{loc}}{(\Lambda_0, \underbar{\bf{0}})}$. We denote $v^+_{\Lambda_0, \underbar{\bf{0}}}$ by $v^+_{\Lambda_0}$. The following relations hold in $W^+_{\mathrm{loc}}{(\Lambda_0, \underbar{\bf{0}})}$:
 
	\begin{center} 
	$e_{i,\underbar{\textbf{k}}}.\hspace{0.05cm} v_{\Lambda_0}^+=0 \,\,\forall  \,\,i \in I,  \forall \,\, \underbar{\textbf{k}} \in \Z^{n-1}_{\geq 0},\,\,\,\,\,\,\,\,\, f_{i,\underbar{\textbf{k}}} v_{\Lambda_0}^+=0 \,\,\forall  \,\,i \in I,  \forall \,\, \underbar{\textbf{k}} \in \Z^{n-1}_{\geq 0}$,\\
	\vspace{0.50cm}
	 $\hspace{1.15cm}h.v_{\Lambda_0}^+= \Lambda_{0}(h)v_{\Lambda_0}^+ \,\, \forall \,\, h \in \fh_{\mathrm{aff}},$ \hspace{0.75cm}
	$h_{i,\underbar{\textbf{k}}}.\hspace{0.05cm} v_{\Lambda_0}^+=0 \,\, \forall i \in I, \,\, \forall \,\,\underbar{\textbf{k}} \in \Z^{n-1}_{\geq 0} $, \\ 
	\vspace{0.25cm}
	$\hspace{-1.25cm} f_{i}. \hspace{0.05cm}v_{\Lambda_0}^+=0 \,\, \forall \,\,i \in I$, \hspace{1.5cm} $f_{0}^2. \hspace{0.05cm}v_{\Lambda_0}^+=0,$
		\vspace{0.25cm}
	$$\hspace{-0.35cm}h_{0, \underbar{\textbf{k}}}\, v^+_{\Lambda_0} = 0\,\, \forall \,\, \underbar{\textbf{k}} \in \mathbb{Z}^{n-1}_{\geq 0} \setminus \{\underbar{\bf{0}}\}, \hspace{0.35cm} f_{0, \underbar{\textbf{k}}}\, v^+_{\Lambda_0} = 0\,\, \forall \,\, \underbar{\textbf{k}} \in \mathbb{Z}^{n-1}_{\geq 0} \setminus \{\underbar{\bf{0}}\}.$$
\end{center}
The relation $f_{0, \underbar{\textbf{k}}}\, v^+_{\Lambda_0} = 0\,\, \forall \,\, \underbar{\textbf{k}} \in \mathbb{Z}^{n-1}_{\geq 0} \setminus \{\underbar{\bf{0}}\}$ follows by calculating the action of $[ e_{0, \underbar{\textbf{k}}}, f_0^2]$ on  $v^+_{\Lambda_0}$ and using
$h_{0, \underbar{\textbf{k}}}\, v^+_{\Lambda_0} = 0\,\, \forall \,\, \underbar{\textbf{k}} \in \mathbb{Z}^{n-1}_{\geq 0} \setminus \{\underbar{\bf{0}}\}$.

\par
	The first step towards an upper bound for $\mathrm{ch}_{q_1}W^{+}_{\mathrm{loc}}{(\Lambda, \underbar{\bf{0}}})$ is the following:

\begin{lemma}\label{action}
Let $\underbar{\bf{a}} = (a_2, a_3, \cdots, a_n) \in \mathbb{Z}_{\geq 0}^{n-1}$ such that not all $a_i$'s zero and let $K=\displaystyle\sum_{i=2}^n a_i$. Then we have

$(\mathrm{i})$
\begin{equation*}
	(e_{\theta} \otimes t_1^{-r}t ^{\underbar{\bf{a}}}) \hspace{0.05cm}v_{\Lambda_0}^+= 
	\begin{cases}
		0& \text{if } r\leq K,\\
		\displaystyle\sum_{m=1}^{r-K} \hspace{0.1cm} (t_1^{-r+m} t^{\underbar{\bf{a}}} K_1) (e_{\theta} \otimes t_1^{-m}) \hspace{0.05cm} v_{\Lambda_0}^+              & \text{if } r>K.\\
	\end{cases}
\end{equation*}

$(\mathrm{ii})$ Further if  $K \geq 2$, then for $j=2,3,\cdots ,n$ we have
\begin{equation*}
	(t_1^{-r} t^{\underbar{\bf{a}}}K_j)\hspace{0.05cm}v_{\Lambda_0}^+= 
	\begin{cases}
		0& \text{if } r\leq K-1,\\
		\displaystyle\sum_{m=1}^{r-(K-1)} (t_1^{-r+m} t_2^{a_2} \cdots t_j^{a_j-1}\cdots t_n^{a_n} K_1) (t_1^{-m}t_j K_j) \hspace{0.05cm} v_{\Lambda_0}^+              & \text{if } r>K-1.\\
	\end{cases}
\end{equation*}
 \end{lemma}
	
\begin{proof}
	We prove both the assertions by induction on $r$. For $r \leq 0$, (i) follows directly from the definition. For $r=1$, (i) follows since $e_{\theta} \otimes t_1^{-1}t_2^{a_2}\cdots t_n^{a_n}= f_{0, \underbar{\textbf{a}}}$ kills $v_{\Lambda_0}^+$. Since $a_j \geq 1$, we have 
	\begin{equation}\label{1s}
	(t_1^{-r}t^{\underbar{\textbf{a}}}K_j) v_{\Lambda_0}^+= \Big([f_{\theta} \otimes t_j,\hspace{0.1cm} e_{\theta} \otimes t_1^{-r}t_2^{a_2}\cdots t_j^{a_j-1} \cdots t_n^{a_n}]-[f_{\theta},\hspace{0.1cm}e_{\theta} \otimes t_1^{-r}t_2^{a_2}\cdots t_j^{a_j} \cdots t_n^{a_n}]\Big) v_{\Lambda_0}^+.
	\end{equation}
	Then (ii) follows for $r\leq 1$ as well. 
	Now let $r\geq 2$ and assume that (i) and (ii) holds for all $r'<r$. Now by Lemma \ref{ab} we have
	\begin{equation}\label{ac}
	\begin{split}
	(e_{\theta} \otimes t_1^{-r} \hspace{0.05cm} t^{\underbar{\textbf{a}}}) \hspace{0.05cm}v_{\Lambda_0}^+ &=\Upsilon_0 \Upsilon_{\theta}\hspace{0.1cm}(e_{\theta} \otimes t_1^{-r+2}\hspace{0.05cm} t^{\underbar{\textbf{a}}}) \hspace{0.05cm}v_{\Lambda_0}^+=\Upsilon_0 \Upsilon_{\theta}\hspace{0.1cm}\Big(e_{\theta} \otimes t_1^{-r+2}\hspace{0.05cm} t^{\underbar{\textbf{a}}} \hspace{0.2cm} \Upsilon_{\theta}^{-1}\Upsilon_{0}^{-1} \hspace{0.05cm}v_{\Lambda_0}^+\Big)\\
	&=\Upsilon_0 \Upsilon_{\theta}\hspace{0.1cm}\Big(e_{\theta} \otimes t_1^{-r+2}\hspace{0.05cm} t^{\underbar{\textbf{a}}} \hspace{0.2cm} \Upsilon_{\theta}^{-1}(f_0v_{\Lambda_0}^+)\Big)=\Upsilon_0 \Upsilon_{\theta}\hspace{0.1cm}\Big(e_{\theta} \otimes t_1^{-r+2}\hspace{0.05cm} t^{\underbar{\textbf{a}}} \hspace{0.2cm} \Upsilon_{\theta}^{-1}(f_0)v_{\Lambda_0}^+\Big)\\
	& = \Upsilon_0 \Upsilon_{\theta}\hspace{0.1cm}\Big((\Upsilon_{\theta}^{-1}(f_0))(e_{\theta} \otimes t_1^{-r+2} \hspace{0.05cm} t^{\underbar{\textbf{a}}})\hspace{0.05cm}v_{\Lambda_0}^+ + [e_{\theta} \otimes t_1^{-r+2}\hspace{0.05cm} t^{\underbar{\textbf{a}}},\hspace{0.05cm} \Upsilon_{\theta}^{-1}(f_0)]\hspace{0.05cm} v_{\Lambda_0}^+\Big) \\
	& = \Upsilon_0 \Upsilon_{\theta}\hspace{0.1cm}\Big((\Upsilon_{\theta}^{-1}(f_0))(e_{\theta} \otimes t_1^{-r+2}\hspace{0.05cm} t^{\underbar{\textbf{a}}})v_{\Lambda_0}^++[e_{\theta} \otimes t_1^{-r+2}\hspace{0.05cm} t^{\underbar{\textbf{a}}}, \hspace{0.05cm}-f_{\theta} \otimes t_1^{-1}] \hspace{0.05cm}v_{\Lambda_0}^+\Big)\\
	& = \Upsilon_0 \Upsilon_{\theta}\hspace{0.1cm}\Big((\Upsilon_{\theta}^{-1}(f_0))(e_{\theta} \otimes t_1^{-r+2}\hspace{0.05cm} t^{\underbar{\textbf{a}}})v_{\Lambda_0}^+\\
	                       &\hspace{4cm}+ [f_{\theta},\hspace{0.05cm} e_{\theta} \otimes t_1^{-r+1}\hspace{0.05cm} t^{\underbar{\textbf{a}}}] \hspace{0.05cm}v_{\Lambda_0}^+ - t_1^{-r+1}\hspace{0.05cm} t^{\underbar{\textbf{a}}} K_1 v_{\Lambda_0}^+\Big)\\
	& = \Upsilon_0 \Upsilon_{\theta}\hspace{0.1cm}\Big( X+Y -t_1^{-r+1} \hspace{0.05cm} t^{\underbar{\textbf{a}}} K_1 v_{\Lambda_0}^+\Big),
\end{split}
	\end{equation}
	where,
	\begin{equation}\label{ad}
	\begin{split}
	X &= \Upsilon_{\theta}^{-1}(f_0)(e_{\theta} \otimes t_1^{-r+2} \hspace{0.05cm} t^{\underbar{\textbf{a}}})\hspace{0.05cm}v_{\Lambda_0}^+\\
	&= \begin{cases}
	\Upsilon_{\theta}^{-1}(f_0) \displaystyle\sum_{m=1}^{r-K-2} \hspace{0.1cm} t_1^{-r+m+2} \hspace{0.05cm} t^{\underbar{\textbf{a}}} K_1 (e_{\theta} \otimes t_1^{-m}) \hspace{0.05cm} v_{\Lambda_0}^+ & \text{if }  r>K+2\\
	0 & \text{if } r \leq K+2,
	\end{cases}
	\end{split}
	\end{equation}
	and \\
	\begin{equation}\label{ae}
	\begin{split}
	Y &= [\hspace{0.05cm}f_{\theta}, \hspace{0.05cm} e_{\theta} \otimes t_1^{-r+1}\hspace{0.05cm} t^{\underbar{\textbf{a}}}]\hspace{0.05cm}\hspace{0.05cm} v_{\Lambda_0}^+= f_{\theta}\hspace{0.05cm}(e_{\theta} \otimes t_1^{-r+1}\hspace{0.05cm} t^{\underbar{\textbf{a}}})\hspace{0.05cm} v_{\Lambda_0}^+\\
	&= \begin{cases}
	f_{\theta} \displaystyle\sum_{m=1}^{r-K-1} \hspace{0.1cm} t_1^{-r+m+1} \hspace{0.05cm}t^{\underbar{\textbf{a}}} K_1 \hspace{0.05cm}(e_{\theta} \otimes t_1^{-m}) \hspace{0.05cm} v_{\Lambda_0}^+ & \text{if }  r>K+1\\
	0 & \text{if } r \leq K+1
	\end{cases}\\
	&= \begin{cases}
	f_{\theta} \displaystyle\sum_{m=0}^{r-K-2} \hspace{0.1cm} t_1^{-r+m+2}\hspace{0.05cm} t^{\underbar{\textbf{a}}} K_1 (e_{\theta} \otimes t_1^{-m-1}) \hspace{0.05cm} v_{\Lambda_0}^+ & \text{if }  r>K+1\\
	0 & \text{if } r \leq K+1.
	\end{cases}
	\end{split}
	\end{equation} The last two expressions follow from the induction hypothesis.
	Now for $r=K+1$, \eqref{ac} becomes
	\begin{equation*}
	\begin{split}
	(e_{\theta} \otimes t_1^{-r} \hspace{0.05cm} t^{\underbar{\textbf{a}}})\hspace{0.05cm} v_{\Lambda_0}^+ & =  \Upsilon_0 \Upsilon_{\theta}\hspace{0.1cm}\Big(-t_1^{-r+1} \hspace{0.05cm} t^{\underbar{\textbf{a}}} K_1 v_{\Lambda_0}^+\Big)\\
	& =-t_1^{-r+1} t^{\underbar{\textbf{a}}} K_1 \Upsilon_0 \Upsilon_{\theta}\hspace{0.05cm}(v_{\Lambda_0}^+)\\
	& =t_1^{-r+1} \hspace{0.05cm}t^{\underbar{\textbf{a}}} K_1 f_0 v_{\Lambda_0}^+\\
	& =t_1^{-r+1} \hspace{0.05cm}t^{\underbar{\textbf{a}}} K_1(e_{\theta}\otimes t_1^{-1})\hspace{0.05cm} v_{\Lambda_0}^+.
	\end{split}
	\end{equation*}
	So (i) holds for $r=K+1$. Now for $r \leq K$, \eqref{ac} becomes
	$$(e_{\theta} \otimes t_1^{-r} \hspace{0.05cm}t^{\underbar{\textbf{a}}})\hspace{0.05cm} v_{\Lambda_0}^+=\Upsilon_0 \Upsilon_{\theta}\hspace{0.1cm}\Big(-t_1^{-r+1} \hspace{0.05cm} t^{\underbar{\textbf{a}}} K_1 v_{\Lambda_0}^+\Big)=0$$ by induction, since $r-1 \leq K-1$ and $t_1^{-r+1} \hspace{0.05cm} t^{\underbar{\textbf{a}}} K_1= \displaystyle\sum_{j=2}^{n} \hspace{0.1cm} \frac{a_j}{r-1} t_1^{-r+1} \hspace{0.05cm} t^{\underbar{\textbf{a}}} K_j$. So (i) holds for $r \leq K$ and hence (i) holds for $r \leq K+1$. Now we consider the case where $r \geq K+2$.
	Note that \eqref{ac} is still valid for $a_2=a_3=\cdots a_n=0$ and $q=m+2$. Hence we have \begin{equation}\label{at}
		(e_{\theta} \otimes t_1^{-m-2}) \hspace{0.05cm}v_{\Lambda_0}^+=\Upsilon_0 \Upsilon_{\theta}\hspace{0.1cm}\Big((\Upsilon_{\theta}^{-1}(f_0))(e_{\theta} \otimes t_1^{-m})v_{\Lambda_0}^++ f_{\theta} \hspace{0.05cm}(e_{\theta} \otimes t_1^{-m-1}) \hspace{0.05cm}v_{\Lambda_0}^+\Big).
	\end{equation}
	In particular for $m=0$, we have 
	\begin{equation}\label{as}
		(e_{\theta} \otimes t_1^{-2})\hspace{0.05cm}v_{\Lambda_0}^+=\Upsilon_0 \Upsilon_{\theta}\hspace{0.1cm}\Big((\Upsilon_{\theta}^{-1}(f_0))(e_{\theta}v_{\Lambda_0}^+)+ f_{\theta}\hspace{0.05cm}(e_{\theta} \otimes t_1^{-1}) \hspace{0.05cm}v_{\Lambda_0}^+)\Big)=\Upsilon_0 \Upsilon_{\theta}\hspace{0.1cm}\Big(f_{\theta} \hspace{0.05cm}(e_{\theta} \otimes t_1^{-1})\hspace{0.05cm} v_{\Lambda_0}^+\Big).
	\end{equation}
	Then by using \eqref{ac}, \eqref{ad}, \eqref{ae}, \eqref{at}, the equation \eqref{as} becomes
	\begin{equation*}
		\begin{split}
	(e_{\theta} \otimes t_1^{-r} \hspace{0.05cm} t^{\underbar{\textbf{a}}}) \hspace{0.05cm}v_{\Lambda_0}^+&= \displaystyle\sum_{m=1}^{r-K-2} \hspace{0.1cm} t_1^{-r+m+2} \hspace{0.05cm}t^{\underbar{\textbf{a}}} K_1 \hspace{0.1cm} \Upsilon_0 \Upsilon_{\theta} \Big(\Upsilon_{\theta}^{-1}(f_0))(e_{\theta} \otimes t_1^{-m})v_{\Lambda_0}^++f_{\theta}\hspace{0.05cm}(e_{\theta} \otimes t_1^{-m-1}) \hspace{0.05cm}v_{\Lambda_0}^+\Big)\\
	&\hspace{2.5cm}+t_1^{-r+2} \hspace{0.05cm}t^{\underbar{\textbf{a}}} K_1 \hspace{0.1cm} \Upsilon_0 \Upsilon_{\theta} \Big(f_{\theta}\hspace{0.05cm}(e_{\theta} \otimes t_1^{-1})\hspace{0.05cm} v_{\Lambda_0}^+\Big)- t_1^{-r+1} \hspace{0.05cm}t^{\underbar{\textbf{a}}} K_1 \hspace{0.1cm} \Upsilon_0 \Upsilon_{\theta}(v_{\Lambda_0}^+)\\
	&=\displaystyle\sum_{m=1}^{r-K} \hspace{0.1cm} t_1^{-r+m}\hspace{0.05cm} t^{\underbar{\textbf{a}}} K_1 (e_{\theta} \otimes t_1^{-m}) \hspace{0.05cm} v_{\Lambda_0}^+.
\end{split}
	\end{equation*}
	So (i) is proved. 
	
	Using \eqref{1s}, we have
	\begin{multline}\label{af}
	(t_1^{-r}t_2^{a_2}\cdots t_j^{a_j}\cdots t_n^{a_n}K_j)\hspace{0.05cm} v_{\Lambda_0}^+  = (f_{\theta} \otimes t_j)\hspace{0.05cm} (e_{\theta} \otimes t_1^{-r}t_2^{a_2}\cdots t_j^{a_j-1} \cdots t_n^{a_n}) \hspace{0.05cm} v_{\Lambda_0}^+\\ -f_{\theta}\hspace{0.05cm}(e_{\theta} \otimes t_1^{-r}t_2^{a_2}\cdots t_j^{a_j} \cdots t_n^{a_n}) \hspace{0.05cm} v_{\Lambda_0}^+.
	\end{multline}
Then using \eqref{af} and by the induction hypothesis, we have $(t_1^{-r}t^{\underbar{\textbf{a}}}K_j)\hspace{0.05cm} v_{\Lambda_0}^+ =0$ if $r < K$ and  for $r= K$ we have 
\begin{equation} 
\begin{split} 
(t_1^{-r}t^{\underbar{\textbf{a}}}K_j)\hspace{0.05cm} v_{\Lambda_0}^+ & =(f_{\theta} \otimes t_j) (t_1^{-K+1}t_2^{a_2}\cdots t_j^{a_j-1} \cdots t_n^{a_n} K_1)(e_{\theta} \otimes t_1^{-1})\hspace{0.05cm}v_{\Lambda_0}^+\\
		& =(t_1^{-K+1}t_2^{a_2}\cdots t_j^{a_j-1} \cdots t_n^{a_n} K_1)\hspace{0.05cm} [f_{\theta} \otimes t_j, e_{\theta} \otimes t_1^{-1}]\hspace{0.05cm}v_{\Lambda_0}^+ \\
		& =(t_1^{-K+1}t_2^{a_2}\cdots t_j^{a_j-1} \cdots t_n^{a_n} K_1) \Big([f_{\theta} , e_{\theta} \otimes t_1^{-1}t_j]+ t_1^{-1}t_j K_j\Big)\hspace{0.05cm}v_{\Lambda_0}^+\\
		& =(t_1^{-K+1}t_2^{a_2}\cdots t_j^{a_j-1} \cdots t_n^{a_n} K_1) ( t_1^{-1}t_j K_j)\hspace{0.05cm}v_{\Lambda_0}^+ .
\end{split}
\end{equation}

So (ii) holds for $q \leq K$. Now we  consider the case where $q>K$. By induction the RHS of \eqref{af} becomes
 \begin{multline}\label{ag}
 	(f_{\theta} \otimes t_j) \displaystyle\sum_{m=1}^{r-(K-1)} \hspace{0.1cm} (t_1^{-r+m}t_2^{a_2}\cdots t_j^{a_j-1} \cdots t_n^{a_n} K_1)(e_{\theta} \otimes t_1^{-m})\hspace{0.05cm}v_{\Lambda_0}^+ \\ 
	- f_{\theta} \displaystyle\sum_{s=1}^{r-K} \hspace{0.1cm} (t_1^{-r+s}t_2^{a_2}\cdots t_j^{a_j} \cdots t_n^{a_n} K_1)(e_{\theta} \otimes t_1^{-s})\hspace{0.05cm}v_{\Lambda_0}^+.
 \end{multline} 
 We have
 \begin{equation}\label{aw}
 \begin{split}
 (f_{\theta} \otimes t_j) (e_{\theta} \otimes t_1^{-m})\hspace{0.05cm}v_{\Lambda_0}^+ & =[f_{\theta} \otimes t_j,\hspace{0.05cm} e_{\theta} \otimes t_1^{-m}]\hspace{0.05cm}v_{\Lambda_0}^+\\ &= \Big([f_{\theta}, \hspace{0.05cm} e_{\theta} \otimes t_1^{-m}t_j]+ t_1^{-m}t_j K_j \Big)\hspace{0.05cm}v_{\Lambda_0}^+\\ & =f_{\theta}\hspace{0.05cm}(e_{\theta} \otimes t_1^{-m}t_j)\hspace{0.05cm} v_{\Lambda_0}^++t_1^{-m}t_j K_j v_{\Lambda_0}^+ \\ & =f_{\theta} \displaystyle\sum_{p=1}^{m-1} \hspace{0.1cm} (t_1^{-m+p}t_j K_1)(e_{\theta} \otimes t_1^{-p})\hspace{0.05cm}v_{\Lambda_0}^++ t_1^{-m}t_j K_j v_{\Lambda_0}^+.
 \end{split}
 \end{equation}
 Therefore the RHS of  \eqref{af} becomes
\begin{multline}\label{ai}
 \displaystyle\sum_{m=1}^{r-(K-1)}  (t_1^{-r+m}t_2^{a_2}\cdots t_j^{a_j-1} \cdots t_n^{a_n} K_1) \hspace{0.05cm}f_{\theta} \displaystyle\sum_{p=1}^{m-1}  (t_1^{-m+p}t_j K_1)(e_{\theta} \otimes t_1^{-p})\hspace{0.05cm}v_{\Lambda_0}^+\\
 + \displaystyle\sum_{m=1}^{r-(K-1)}  \hspace{-0.35cm} (t_1^{-r+m}t_2^{a_2}\cdots t_j^{a_j-1} \cdots t_n^{a_n} K_1) (t_1^{-m}t_j K_j)\hspace{0.05cm} v_{\Lambda_0}^+\\
 - f_{\theta} \displaystyle\sum_{s=1}^{r-K}  \hspace{-0.05cm} (t_1^{-r+s}t_2^{a_2}\cdots t_j^{a_j} \cdots t_n^{a_n} K_1)(e_{\theta} \otimes t_1^{-s}) \hspace{0.05cm}v_{\Lambda_0}^+.
 \end{multline}
By applying $h_{\theta} \otimes t_j$ to both sides of (i), we see that the first and the third terms of \eqref{ai} are equal. Hence (ii) is proved.
\end{proof}

Let $\bar{\mathcal{Z}}$ be the subalgebra of $\mathcal{U}(\mathcal{T}^{+})$ generated by $\{t_1^{-r}t_2^{a_2}\cdots t_j^{a_j}\cdots t_n^{a_n}K_j: r>0, \underbar{\textbf{a}} \in  \mathbb{Z}^{n-1}_{\geq 0}, j=2,3, \cdots, n\}$ and let $\hat{\mathcal{Z}}$ to be the subalgebra generated by the elements $\{t_1^{-r}t_j K_{1}:r>0, j=2,3, \cdots, n\}$.
	\begin{prop}\label{ah}
		We have  $\bar{\mathcal{Z}}v_{\Lambda_0}^+= \hat{\mathcal{Z}}v_{\Lambda_0}^+$.
		\begin{proof} We use induction on $K=\sum_{i=2}^n a_i$ \, to show that the elements $t_1^{-r}t_2^{a_2}\cdots t_j^{a_j}\cdots t_n^{a_n}K_j v_{\Lambda_0}^+$ are in $\hat{\mathcal{Z}}v_{\Lambda_0}^+$. 
			For $K=1$, the elements of $\bar{\mathcal{Z}}$ are of the form $t_1^{-r}t_j K_{j}$ for some $j=2,\cdots, n$. Then the assertion follows as $t_1^{-r}t_j K_{j}=rt_1^{-r}t_j K_{1}$. For $K=2$, the elements of $\bar{\mathcal{Z}}$ are of the form $t_1^{-r}t_j^2 K_{j}$ or $t_1^{-r}t_i t_j K_{j}$ $(i\neq j)$. From Lemma \ref{action}(ii) we have
			\[
			(t_1^{-r}t_j^{2}K_j)\hspace{0.05cm}v_{\Lambda_0}^+= 
			\begin{cases}
				0& \text{if } r\leq 1,\\
				\displaystyle\sum_{s=1}^{r-1} (t_1^{-r+s} t_j K_1) (t_1^{-s}t_j K_j) \hspace{0.05cm} v_{\Lambda_0}^+              & \text{if } r>1\\
			\end{cases}
			\]
			and
			\[
			(t_1^{-r}t_i t_jK_j)\hspace{0.05cm}v_{\Lambda_0}^+= 
			\begin{cases}
			0& \text{if } r\leq 1,\\
			\displaystyle\sum_{s=1}^{r-1} (t_1^{-r+s} t_i K_1) (t_1^{-s}t_j K_j) \hspace{0.05cm} v_{\Lambda_0}^+              & \text{if } r>1.\\
			\end{cases}
			\]
			So for $K=2$ the result follows. Now let $K> 2$ and assume that the result holds for any natural number less than $K$. Consider $t_1^{-r}t_2^{a_2}\cdots t_j^{a_j}\cdots t_n^{a_n}K_j v_{\Lambda_0}^+$ with $\displaystyle \sum_{i=2}^n a_i=K$. From Lemma \ref{action}(ii), we have
			\[
			(t_1^{-r}t^{\underbar{\textbf{a}} }K_j)\hspace{0.05cm}v_{\Lambda_0}^+= 
			\begin{cases}
			0& \text{if } r\leq K-1,\\
			\displaystyle\sum_{s=1}^{r-(K-1)} (t_1^{-r+s} t_2^{a_2} \cdots t_j^{a_j-1}\cdots t_n^{a_n} K_1) (t_1^{-s}t_j K_j) \hspace{0.05cm} v_{\Lambda_0}^+              & \text{if } r>K-1.\\
			\end{cases}
			\]
			Since 
			$$t_1^{-r+s} t_2^{a_2} \cdots t_j^{a_j-1}\cdots t_n^{a_n} K_1=
			\displaystyle\sum_{i=2, i\neq j}^{n} \frac{a_i}{r-s} t_1^{-r+s} t_2^{a_2} \cdots t_j^{a_j-1}\cdots t_n^{a_n} K_i $$ $$ +\,\, \frac{a_{j}-1}{r-s} t_1^{-r+s} t_2^{a_2} \cdots t_j^{a_j-1}\cdots t_n^{a_n} K_j,$$
		the assertion holds using induction. 
			\end{proof}
	\end{prop}
	
\begin{prop} \label{pmp}
	We have $W_{\mathrm{loc}}(\Lambda_0)= \hat{\mathcal{Z}} \, \mathcal{U}(n_{\mathrm{aff}}^{-})v_{\Lambda_0}^+$.
	\begin{proof}
		Combining Lemma \ref{action}  and Proposition \ref{ah} the proof follows.
		\end{proof}
	\end{prop}

	Using the above results, for ${\underbar{\bf{a}}} \in ({\mathbb{C}^{*})}^{n-1}$ we have the following proposition. 

\begin{prop}\label{main}
		$$\mathrm{ch}_{q_1} W_{\mathrm{loc}}(\Lambda_0, {\underbar{\bf{a}}}) \leq  \mathrm{ch}_{q_1} W^+_{\mathrm{loc}}(\Lambda_0, {\underbar{\bf{a}}}) \leq \hspace{0.2cm} \mathrm{ch}_{q_1} L(\Lambda_0)\Big(\prod_{m>0}{\frac{1}{1- {q_1}^m}}\Big)^{n-1},$$
	$$\mathrm{ch}_{q_1, q_2, \ldots, q_{n}} W_{\mathrm{loc}}(\Lambda_0, {\underbar{\bf{a}}}) \leq \mathrm{ch}_{q_1, q_2, \ldots, q_{n}} W^+_{\mathrm{loc}}(\Lambda_0, {\underbar{\bf{a}}})  \leq \hspace{0.2cm} {\displaystyle{ \mathrm
			{ch}_{q_1} L(\Lambda_0)}\prod_{m>0, i = 2}^{n}{\frac{1}{1- {q_1}^m q_i}}}.$$
\end{prop}

\begin{proof}
For ${\underbar{\bf{a}}} \in ({\mathbb{C}^{*})}^{n-1}$, consider $ W_{\mathrm{loc}}(\Lambda_0, {\underbar{\bf{a}}})$ as a $\CT^{+}$-module by a restriction. Then by Corollary \ref{c1} the map $W^+_{\mathrm{loc}}(\Lambda_0, {\underbar{\bf{a}}}) \rightarrow
W_{\mathrm{loc}}(\Lambda_0, {\underbar{\bf{a}}})$ given by $v^+_{\Lambda_0, \underbar{\bf{a}}} \mapsto v_{\Lambda_0, \underbar{\bf{a}}}$
is a surjective $\CT^{+}$-module map.  Then the first inequalities of characters of the above two expressions follow. Proposition \ref{pmp} settles the second inequality.
\end{proof}

In the next section we prove the equalities of the characters by using representation theory of Fock space of $\CT$.

\section{Fock space and vertex operators} \label{s4}
In this section we take our ``affine variable" as ``$\,t_{n}$" instead of ``$\,t_{1}$" and set $\mathcal{T}_{\mathrm{aff}}^{n}=\mathfrak{g}\otimes \C[t_{n}^{\pm 1}] \oplus \C K_n \oplus \C d_n$.  We use notations from \cite{CMP}. Let $\Gamma$ and $Q$ be root lattices given by \begin{center}
	
$\Gamma= \bigoplus\limits_{i=1}^l \Z \alpha_i  \oplus \bigoplus\limits_{i=1}^{n-1} \Z \delta_i \oplus \bigoplus\limits_{i=1}^{n-1} \Z d_i$ \,\,\,\,\,\, and \,\,\,\,\,\, $Q= \bigoplus\limits_{i=1}^l \Z \alpha_i \oplus \bigoplus\limits_{i=1}^{n-1} \Z \delta_i$.
\end{center}  
The lattice $\Gamma$ possees a non-degenerate form $(,)$ while $Q$ does not (see \cite[section 3]{CMP} for more details). Let $\mathfrak{p}= \C \otimes_{\Z} \Gamma$ and $\mathfrak{h}= \C \otimes_{\Z} Q$, and consider their corresponding Heisenberg algebras 
\begin{center}
 $\hat{\bb}=\bigoplus\limits_{k \in \Z} \mathfrak{p}(k) \oplus \C K_n$ \,\,\, and \,\,\, $\hat{\e}=\bigoplus\limits_{k \in \Z} \mathfrak{h}(k) \oplus \C K_n$,
\end{center}
 where for each $k \in \mathbb Z$, $\mathfrak{p}(k)$ and $\mathfrak{h}(k)$ denote isomorphic copies of $\mathfrak{p}$ and $\mathfrak{h}$ respectively and $K_n$ is the central element of the affine Lie algebra $\mathcal{T}_{\mathrm{aff}}^{n}$. The bracket on  $\hat{\bb}$ (similarly on $\hat{\e}$) is given by 
$$[\alpha(k), \beta(s)]=k(\alpha, \beta)\,  \delta_{k, -s} K_n \hspace{1cm} \mathrm{and}  \hspace{1cm}  [\alpha(k), K_n]=0, \,\,\,\, \forall k \in \mathbb Z.$$

 Set $$\bb=\bigoplus\limits_{\substack{k \in \Z \\ k\neq 0}} \mathfrak{p}(k) \oplus \C K_n,  \,\,\, \e=\bigoplus\limits_{\substack{k \in \Z\\ k\neq 0}} \mathfrak{h}(k) \oplus \C K_n, \,\,\, \bb_{\pm}= \bigoplus\limits_{k \gtrless 0} \mathfrak{p}(k), \,\,\, \e_{\pm}= \bigoplus\limits_{k \gtrless 0} \mathfrak{h}(k).$$\\
Let $S(\bb_{-})$ denote the symmetric algebra of $\bb_{-}$ with the following action of $\bb$:

$K_n$ acts as 1$;\,\, \alpha(-m)$ acts as left multiplication for $m>0$;
$\alpha(m)$ acts as a derivation given by $\alpha(m) (\beta(-n))= m \delta_{m,n} (\alpha, \beta)$ for $\alpha, \beta \in \Gamma, m,n>0$. 
With the above action $S(\bb_{-})$ is an irreducible representation of $\bb$.

 For $\um=(m_1,m_2,\cdots m_{n-1})$, we set $\delta_{\um}=m_1 \delta_1+ m_2 \delta_2+ \cdots, m_{n-1} \delta_{n-1}$. Then we have $(\delta_{\um}, \delta_{\um})=0$. 
  
 Let $\epsilon$ be a 2-cocycle on $Q_{fin}$ as defined in \cite[(2.11)]{CMP}. Then $\epsilon$ can be extended to $Q$ by defining $\epsilon(\alpha, \delta_{\um})=1$ for all $\alpha \in Q,\,  \um \in \Z^{n-1}$. Further $\epsilon$ is extended from $Q \times Q$ to $Q \times \Gamma$ in a appropriate way. Let $\C(\Gamma)$ and $\C(Q)$ denote the group algebras of $\Gamma$ and $Q$ respectively. For $x\in \Gamma$, $e^x$ denotes its image in $\C(\Gamma)$. The group algebra $\C(\Gamma)$ is a $\C(Q)$ module with the action $e^{\alpha} e^{\gamma}= \epsilon (\alpha, \gamma) e^{\alpha+ \gamma}$ for $\alpha \in Q, \gamma \in \Gamma$.\par \vspace{0.2cm}
 
 For any $\e_{-}$ submodule $M \subseteq S(\bb_{-})$, $\C(\Gamma) \otimes M$ is a $\e$ module with the 
 action $a(m) (e^{\gamma} \otimes v)= e^{\gamma} \otimes a(m)v$ for $a \in \mathfrak{h}$ and $m\neq 0$ and it is extended to $\hat{\e}$ by $a(0)\, (e^{\gamma} \otimes v)= (a,\gamma)\,\, e^{\gamma} \otimes v$,  $a \in \mathfrak{h}$.\\
 Let us take $M=S(\e_{-})$ and consider the $\hat{\e}$ module $V= \C(\Gamma)\otimes S(\e_{-})$.
 To extend the action of $\hat{\e}$ on $V$ to $\mathcal{T}$, we need to recall the vertex operators. \par
  Let $\alpha \in Q$ and $z$ be a variable, define $T_{\pm}(\alpha, z) := - \sum_{r \gtrless 0} \frac{\alpha(r) }{r} z^{-r}$.
 	The vertex operator for $\alpha \in Q$ is defined as
 	$$ X(\alpha, z)= z^{\frac{(\alpha, \alpha)}{2}} \,\, \text{exp}\,\, T_{-}(\alpha, z) \,\, e^{\alpha} \,\, z^{\alpha(0)} \,\,\text{exp} \,\, T_{+}(\alpha, z).$$
 	The operator $z^{\alpha(0)}$ is defined on $V$ by $z^{\alpha(0)} (e^{\gamma} \otimes v)= z^{(\alpha, \gamma)}\,  e^{\gamma} \otimes v$. Let $X_r(\alpha) $ be the coefficient of $z^{-r}$ in the expression of $X(\alpha, z)$.  Then we can write $X(\alpha, z)$ in the form $X(\alpha, z)= \displaystyle\sum_{r \in \Z} X_r(\alpha) z^{-r}$.
 	For $\alpha \in Q$, define\\ 
 	$$\alpha^{+}(z)= \displaystyle\sum_{r  \geq 0} \alpha(r) z^{-r-1},\,\,\,\, \alpha^{-}(z)= \displaystyle\sum_{r < 0} \alpha(r) z^{-r-1}, \,\,\,\, \alpha(z)= \alpha^{+}(z) + \alpha^{-}(z) =\displaystyle\sum_{r \in \Z} \alpha(r) z^{-r-1}.$$
 	For $\alpha \in Q$ and $\delta_{\um}$, define  $$T_{\delta_{\um}, z}^{\alpha}=\,\, :\alpha(z) X({\delta_{\um}, z}):\,\, = \alpha^{-}(z) X({\delta_{\um}, z}) + X({\delta_{\um}, z})\, \alpha^{+}(z).$$ 
 	We define $T_{r}^{\alpha} (\delta_{\um})$ by the equation $T_{\delta_{\um}, z}^{\alpha}= \displaystyle\sum_{r \in \Z} T_{r}^{\alpha} (\delta_{\um}) z^{-r-1}$. \par
 	 Now we state a result which extends the action of $\hat{\e}$ to $\mathcal{T}$ on $V:= \C(\Gamma) \otimes  S(\e_{-})$.
 	 
 \begin{theorem}[\cite{CMP}, Theorem 3.14]
 	 The following map extend the action of $\hat{\e}$ to  $\mathcal{T}$ on $V= \C(\Gamma) \otimes  S(\e_{-})$. 
\begin{equation*}
	\begin{split}
	X_{\alpha} \otimes t^{\textbf{m}} & \mapsto X_{m_n}(\alpha+ \delta_{\um}) \hspace{0.5cm} \alpha \in R_{fin} \\
		h \otimes t^{\textbf{m}} & \mapsto T_{m_n}^{h}(\delta_{\um}) \hspace{1.35cm} h \in \mathfrak{h}\\
		t^{\textbf{m}} K_i & \mapsto T_{m_n}^{\delta_{i}}(\delta_{\um}) \hspace{1.35cm} 1 \leq i \leq n-1 \\
		t^{\textbf{m}} K_n & \mapsto X_{m_n}(\delta_{\um})
	\end{split}
\end{equation*}
The action of $d_i$'s are given by 
\begin{equation*}
	\begin{split}
	\hspace{1cm}d_i\, (e^{\gamma} \otimes \alpha(r))	 & = (d_i, \gamma) (e^{\gamma} \otimes \alpha(r)),\hspace{1.5cm} 1 \leq i \leq n-1 \\
		d_n(e^{\gamma} \otimes \alpha(r)) & =(\frac{(\gamma, \gamma)}{2}+r) (e^{\gamma} \otimes \alpha(r)).
	\end{split}
\end{equation*}
 	\end{theorem}
 \vspace{0.3cm}
 \subsection{The submodule ${\bf{V(0)}}$}
We have $V= \C(\Gamma) \otimes  S(\e_{-})=\oplus_{\gamma \in \Gamma/ Q} \, V(\gamma)$, where $V(\gamma)= e^{\gamma+Q} \otimes S(\e_{-})$ is a cyclic $\mathcal{T}$ module with generator $e^{\gamma} \otimes 1$. Our main module of study is $$V(0)= e^{Q} \otimes S(\e_{-}).$$ Set 
$$ \overset{\circ}{\e}_{-}= \bigoplus\limits_{k < 0} {\mathfrak{h}}(k).$$ 
By \cite[section 5]{1995} and \cite{FK} we know that $\C[Q_{fin}] \otimes\, S(\overset{\circ}{\e}_{-})$ is isomorphic to the level one integrable irreducible module $L^n(\Lambda_{0})$. In the next theorem, we describe the structure of $V(0)$. 

\begin{prop}
	As a vector space $V(0) \cong \C[Q_{fin}] \otimes S(\overset{\circ}{\e}_{-}) \otimes M \otimes \C[\tau_{1}^{\pm 1}, \tau_{2}^{\pm 1} ,\cdots, \tau_{n-1}^{\pm 1}]$, where $M$ is the polynomial algebra generated by $\{\delta_{1}(k), \delta_{2}(k), \cdots, \delta_{n-1}(k): k<0\}$ and $\tau_i= X_0(\delta_i)$.
	\end{prop}
\begin{proof}
	Follows from the Lemmas 5.1, 5.2, 5.3 of \cite{1995}.
	\end{proof}
	
	As a corollary we get the following;
\begin{cor}\label{bf}
	The module $V(0)$ is free over $M \otimes \C[\tau_{1}^{\pm 1}, \tau_{2}^{\pm 1} ,\cdots, \tau_{n-1}^{\pm 1}]$. In particular, $V(0)$ is free over $\C[\tau_{1}^{\pm 1}, \tau_{2}^{\pm 1} ,\cdots, \tau_{n-1}^{\pm 1}]$.
\end{cor}
From above we get that $V(0) \cong \C[Q_{fin}] \otimes S(\overset{\circ}{\e}_{-}) \otimes M \otimes \C[\tau_{1}^{\pm 1}, \tau_{2}^{\pm 1} ,\cdots, \tau_{n-1}^{\pm 1}] \cong L^n(\Lambda_{0}) \otimes M \otimes \C[\tau_{1}^{\pm 1}, \tau_{2}^{\pm 1} ,\cdots, \tau_{n-1}^{\pm 1}]$.
 Now by \cite[Proposition 4.3]{yokunama}, we obtain that the operators $\{X_n(\alpha): \alpha \in R_{fin}, n \in \Z\}$ act on $L^n(\Lambda_{0}) \cong \C[Q_{fin}] \otimes S(\overset{\circ}{\e}_{-})$ by the correspondence
\begin{equation}\label{az}
	e_i \otimes t_n^{k} \mapsto X_k(\alpha), \hspace{2cm} f_i \otimes t_n^{k} \mapsto X_k(-\alpha).
\end{equation}

   Note that $V(0)$ is generated by $v_n \otimes 1 \otimes 1$, where $v_n$ is the highest weight vector for $L^n(\Lambda_{0})$.  In the next theorem we give the $\mathcal{T}$ action on $V(0)$. For $ x\in \mathfrak{g}$ define $x(z)= \displaystyle\sum_{k \in \Z} (x \otimes t_n^k) z^{-k}$, and for $ \um \in \Z^{n-1}$ define $\Delta_{\um}(z)= \prod_{i = 1}^{n-1} \mathrm{exp}\,\,\Big(\displaystyle\sum_{k >0} \frac{m_i \delta_i(-k)}{k} z^k\Big)$ and $\tau^{\um}=\tau_{1}^{m_1}\tau_{2}^{m_2} ,\cdots \tau_{n-1}^{m_{n-1}}$. Let $d^{(i)}$ be the operator which counts the degree of $\tau_{i}$ on $\C[\tau_{1}^{\pm 1}, \tau_{2}^{\pm 1} ,\cdots, \tau_{n-1}^{\pm 1}]$. We make $M$ a graded algebra by defining deg $\delta_i(k)=k$ for $1\leq i\leq n-1$. Let $d^{M}$ be the operator that counts the degree on $M$. The following result is a consequence of Lemma 5.7 of \cite{1995}.
   
\begin{theorem}\label{vr} 
	$V(0)$ is a $\mathcal{T}$ module with the action given by
	\begin{equation*}
		\begin{split}
			\displaystyle\sum_{s \in \Z} x \otimes t^{\um} t_n^s z^{-s} & \mapsto x(z) \otimes \Delta_{\um}(z) \otimes \tau^{\um} \\
			\displaystyle\sum_{s \in \Z} t^{\um} t_n^s K_i z^{-s-1} & \mapsto \text{id} \otimes \delta_i(z) \Delta_{\um}(z) \otimes \tau^{\um}\,\,\,\, ;\,\,\,\,\,\,\, 1\leq i \leq n-1 \\
			\displaystyle\sum_{s \in \Z} t^{\um} t_n^s K_n z^{-s} & \mapsto \text{id} \otimes \Delta_{\um}(z) \otimes \tau^{\um} \\
	\end{split}
	\end{equation*}
\begin{equation*}
	\hspace{-0.5cm}d_i \mapsto \text{id} \otimes \text{id} \otimes d^{(i)}\,\, (1\leq i \leq n-1),
	\hspace{0.5cm}d_n  \mapsto d_n \otimes \text{id} \otimes \text{id} + \text{id} \otimes d^{M} \otimes \text{id}.
\end{equation*}
\end{theorem}
As in \eqref{aut}  we define an automorphism of $\CT$ associated with the following $n \times n$ matrix 
\begin{equation*}
	A = \begin{pmatrix}
		0 & 0 & 0 & \cdots & 0 & -1 \\
		0 & 1 & 0 & \cdots & 0 & 0 \\
		0 & 0 & 1 & \cdots & 0 & 0 \\
		\vdots & \vdots & \vdots  & \ddots & \vdots & \vdots  \\
		0 & 0 & 0 & \cdots & 1 & 0\\
		1 & 0 & 0 & \cdots & 0 &0
	\end{pmatrix} 
=
\begin{pmatrix}
	0 & \rvline &\bigzero_{1 , {n-2}} & \rvline &-1 \\
	\hline 
	\bigzero_{n-2 ,1} & \rvline &\bigi_{{n-2}, {n-2}} & \rvline & \bigzero_{n-2 ,1}\\
	
	\hline
	1 & \rvline &\bigzero_{1 , {n-2}} & \rvline &0
\end{pmatrix}
\end{equation*} 
 which is again denoted by $A$. The image is given by the following
\begin{equation*}
	A(x\otimes t^{\m})= x \otimes t_1^{-m_n} t_2^{m_2}\cdots t_{n-1}^{m_{n-1}} t_n^{m_1},\,\,\, x \in \mathfrak{g}
\end{equation*} 
\begin{center}
	$A(t^{\m} K_1)=  t_1^{-m_n} t_2^{m_2}\cdots t_{n-1}^{m_{n-1}} t_n^{m_1} K_n$,\ 
\end{center}
\begin{equation*}
	A(t^{\m} K_i)=  t_1^{-m_n} t_2^{m_2}\cdots t_{n-1}^{m_{n-1}} t_n^{m_1} K_i , \,\,\,\,\,  2 \leq i\leq n-1
\end{equation*}
\begin{center}
	$A(t^{\m} K_n)=  -t_1^{-m_n} t_2^{m_2}\cdots t_{n-1}^{m_{n-1}} t_n^{m_1} K_1$
\end{center}
\begin{equation*}
	A(d_1)= d_n,\  A(d_i)= d_i,\ (2 \leq i\leq n-1),\ A(d_n)= -d_1.
\end{equation*}
Let $\V= A^* V(0)$, the pull back of $V(0)$ via the automorphism $A$. We denote $\bf{v}$ by the image of $v_n \otimes 1 \otimes 1$.
\begin{lemma}\label{bg}
	We have
	$h_{i, \lk}.\, {\bf{v}}=0$,  $\forall i \in I, \,\, \underbar{\bf{k}} \in \Z^{n-1}$.
	\begin{proof}
		For all $i \in I$, we have
		\begin{equation*}
			h_{i,\underbar{\bf{k}}}\,.{\bf{v}} = A(h_{i,\underbar{\textbf{k}}})(v_n \otimes 1 \otimes 1) = h_i \otimes t_1^{-k_n} t_2^{k_2}\cdots t_{n-1}^{k_n-1} (v_n \otimes 1 \otimes 1).
		\end{equation*}
		Let $\ukk = (-k_{n}, k_2,\cdots, k_{n-1})$ for $\lk = (k_2, k_3, \cdots, k_n)$ and let $\Delta_{\uk}(z)= \sum_{s \geq 0} (\Delta_{\uk})^{(s)} \, z^s$.
		Now since 
		$\displaystyle\sum_{s \in \Z} h_i \otimes t_1^{-k_n} t_2^{k_2}\cdots t_{n-1}^{k_{n-1}} t_n^{s} z^{-s}$ acts on $V(0)$ by $h_i(z) \otimes \Delta_{\ukk}\, (z) \otimes \tau^{\ukk}$, the action of $h_{i,\underbar{\bf{k}}}$ on $V(0)$ is given by 
		$ \sum_{s\geq 0} (h_i \otimes t_n^ s) \otimes (\Delta_{\ukk})^{(s)} \otimes \tau^{\ukk}$. Again since $(h_i \otimes t_n^ s). v_n=0$ for $s\geq 0$, the lemma follows.
		\end{proof}
	\end{lemma}
\begin{theorem}\label{isom}
	The assignment $v_{\Lambda_0} \mapsto \bf{v}$ gives a surjective homomorphism between $W_{\mathrm{glob}}{(\Lambda_0)}$ and $\V$.
\end{theorem} 
\begin{proof}
	The surjectivity of the map follows from the definition of $\V$. So we only need to check that the map is well-defined. It is easy to see that $\mathfrak{g}.(v_n \otimes 1 \otimes 1)= \mathfrak{g}.v_n \otimes 1 \otimes 1$. So $e_i. {\bf{v}}=0=f_i.{\bf{v}}$. Then for $1 \leq i \leq l$, we have $e_{i,\underbar{\bf{k}}}{\bf{v}}=0$ by using the fact that $e_i.{\bf{v}}=0$ and $h_{i,\underbar{\bf{k}}}. {\bf{v}}=0$.
	 Now since $A(e_{0,\underbar{\bf{k}}})= f_{\theta} \otimes t^{\ukk} t_n$, the action of $e_{0,\underbar{\bf{k}}}$ is given by the expression $ \sum_{s=-1}^{\infty} (f_\theta \otimes t_n^ s) \otimes (\Delta_{\ukk})^{(s+1)} \otimes \tau^{\ukk}$. We have \begin{equation*}\label{bh}
	\begin{aligned}
		&\Big((f_\theta \otimes t_n^ {-1}) \otimes (\Delta_{\ukk})^{(0)} \otimes \tau^{\ukk}\Big) (v_n \otimes 1 \otimes 1)\\
		=&(f_\theta \otimes t_n^ {-1}).(v_n) \otimes 1 \otimes \tau^{\ukk}\\
		=& \tau^{\ukk}(f_\theta \otimes t_n^ {-1}).(v_n) \otimes 1 \otimes 1 \,\,\,\,\,\,\, (\text{By Corollary}\, \ref{bf}) \\
		=& \tau^{\ukk} (f_\theta \otimes t_n^ {-1}).(v_n\otimes 1 \otimes 1)\\
		=& 0 \,\,\,\,\,\, (\text{By Lemma}\, \ref{bg}).
	\end{aligned}
\end{equation*}
So $e_{0,\underbar{\bf{k}}}. {\bf{v}}=0$.
  Now $t^{\lm} K_i. {\bf{v}}=0$ ($\lm \in \Z^{n-1}$, $2 \leq i \leq n$) since $(\delta_{i}, \delta_{j})=0$. It is easy to check that $f_{0}^2. {\bf{v}}= f_{0}^2. v_n \otimes 1 \otimes 1=0$.
Also we have 
$$ K_1. {\bf{v}}= K_n(v_n \otimes 1 \otimes 1)= K_n. v_n \otimes 1 \otimes 1= v_n \otimes 1 \otimes 1={\bf{v}}.$$
From the given action, it follows that $d_i. {\bf{v}}=0$ for all $1 \leq i\leq n$. This completes the proof the theorem.
	\end{proof}
 Note that the action of $t_n K_1$ on $\V$ is by $\tau_1^{-1}$ because $A(t_n K_1)= t_1^{-1} K_n$. Similarly the action of $t_i K_1$ corresponds to $\tau_{i}$ for $2 \leq i \leq n-1$. So we can consider $\V$ as an $A(\Lambda_0)= \mathbb{C}[x_1^{\pm1}, \ldots, x_{n-1}^{\pm1}]$ module with the action given by $x_1 \mapsto t_n K_1$, $x_i \mapsto t_i K_1$ for $2 \leq i \leq n-1$. Then $\V$ becomes a free $A(\Lambda_0)$ module by Corollary \ref{bf}. For $\overline{\textbf{a}} \in {(\C^{*})}^{n-1}$, we set $\V_{\overline{\textbf{a}}}:= \V \otimes _{A(\Lambda_0)} \C_ {\overline{\textbf{a}}}.$ We then have the following proposition. 
 
 \begin{prop}
	$\V_{\overline{\textbf{a}}}$ is a  $\bar{\mathcal{T}}$ module with $q_1$-character $$\mathrm{ch}_{q_1}\, \V_{\overline{\textbf{a}}}= \mathrm{ch}_{q_1} L(\Lambda_0)\Big(\prod_{m>0}{\frac{1}{1- {q_1}^m}}\Big)^{n-1}.$$
	\begin{proof}
		This follows since $V(0) \cong L^n(\Lambda_{0}) \otimes M \otimes \C[\tau_{1}^{\pm 1}, \tau_{2}^{\pm 1} ,\cdots, \tau_{n-1}^{\pm 1}]$ and $A(d_1)=d_n$.
		\end{proof}
\end{prop}
\begin{theorem} \label {mainprop}
	 $W_{\mathrm{loc}}(\Lambda_0, {\overline{\textbf{a}}}) \cong \V_{\overline{\textbf{a}}}$ as a $\bar{\mathcal{T}}$-module for every $\overline{\textbf{a}} \in {(\C^{*})}^{n-1}$ and
	\begin{center}
		$\mathrm{ch}_{q_1} W_{\mathrm{loc}}(\Lambda_0, \overline{\textbf{a}}) = \mathrm{ch}_{q_1} W_{\mathrm{loc}}^{+}(\Lambda_0, \overline{\textbf{a}}) =  \mathrm{ch}_{q_1} L(\Lambda_0)\Big(\prod_{m>0}{\frac{1}{1- {q_1}^m}}\Big)^{n-1}.$
	\end{center}

\begin{proof}
	
	By Theorem \ref{isom}, there is a surjective  $\mathcal{T}$-module homomorphism $W_{\mathrm{loc}}(\Lambda_0, {\overline{\textbf{a}}}) \rightarrow \V_{\overline{\textbf{a}}}$ for every $\overline{\textbf{a}} \in {(\C^{*})}^{n-1}$. So we have 
	\vspace{-0.1cm}
\begin{equation}\label{1}
	\text{ch}_{q_1} W_{\mathrm{loc}}^{+}(\Lambda_0, \overline{\textbf{a}}) \geq \text{ch}_{q_1} W_{\mathrm{loc}}(\Lambda_0, \overline{\textbf{a}}) \geq \text{ch}_{q_1}\, \V_{\overline{\textbf{a}}}= \text{ch}_{q_1} L(\Lambda_0)\Big(\prod_{m>0}{\frac{1}{1- {q_1}^m}}\Big)^{n-1}.
\end{equation}
We know that $q_1$-character $\text{ch}_{q_1} W_{\mathrm{loc}}^{+}(\Lambda_0, \overline{\textbf{a}})$ is independent of $\overline{\textbf{a}} \in \C^{n-1}$. So from proposition \ref{main} we have 
\begin{equation}\label{2}
	\text{ch}_{q_1} W_{\mathrm{loc}}^{+}(\Lambda_0, \overline{\textbf{a}}) \leq \text{ch}_{q_1} L(\Lambda_0)\Big(\prod_{m>0}{\frac{1}{1- {q_1}^m}}\Big)^{n-1}=\text{ch}_{q_1}\, \V_{\overline{\textbf{a}}}.
\end{equation} Combining (\ref{1}) and (\ref{2}) we get the desired result.
	\end{proof}
\end{theorem}
\begin{cor}
As a $\mathcal{T}$ module $W_{\mathrm{glob}}(\Lambda_0)$ is isomorphic to  $\V$.
\end{cor}
\begin{proof}
The arguments of the proof are similar to that of Theorem 4.10 of \cite{ko}.
\end{proof}
\begin{cor}
For $\overline{\textbf{a}} \in {(\C^{*})}^{n-1}$ the graded character of $W_{\mathrm{loc}}(\Lambda_0, \overline{\textbf{a}})$ is given by
	$$\mathrm{ch}_{q_1, q_2, \ldots, q_{n}} W_{\mathrm{loc}}(\Lambda_0, \overline{\textbf{a}}) = \hspace{0.2cm} {\displaystyle{ \mathrm{ch}_{q_1} L(\Lambda_0)}\prod_{m>0, i = 2}^{n}{\frac{1}{1- {q_1}^m q_i}}}.$$
	\begin{proof} Follows from Proposition \ref{mainprop} and Theorem \ref{main}.
		
\end{proof}
\end{cor}
\begin{cor}
Let $\CT$ be a toroidal Lie algebra with underlying finite dimensional simple Lie algebra $\mathfrak{g}$ is of the type $A, D$ or $E$. Then
$$\hspace{0.6cm} \mathrm{ch}_{q_1} W_{\mathrm{loc}}(\Lambda_0, \overline{\textbf{a}}) = e^{\Lambda_0} \Big(\prod_{m>0}{\frac{1}{1- {q_1}^m}}\Big)^{n+l-1} \hspace{1cm} \mathrm{and}$$
$$\mathrm{ch}_{q_1, q_2, \ldots, q_{n}} W_{\mathrm{loc}}(\Lambda_0, \overline{\textbf{a}}) = \hspace{0.2cm} e^{\Lambda_0}{\displaystyle{ \prod_{k>0}{\frac{1}{(1-q_1^{k})^l}}}\prod_{m>0, i = 2}^{n}{\frac{1}{1- {q_1}^m q_i}}}.$$
\end{cor}
\begin{proof}
Follows immediately from Proposition 12.13 of \cite{kac}. 
\end{proof}

\section{Appendix on Symmetric functions}

The aim of this appendix is to review some of the basic notions about symmetric functions which we require in the proof of Proposition {\ref{imp}}. Most of the material covered here can be found in \cite{MAC}. Let $x_1, x_2, \cdots, x_n$ be a finite set of indeterminates. The symmetric group on
$n$ letters $S_n$ acts on $\mathbb Z[x_1, x_2, \cdots, x_n]$, the ring of polynomials in $n$ variables by permuting the variables. The ring of symmetric polynomials is defined by 
$ \Lambda(n) := \{ f \in \mathbb{Z}[x_1, \cdots, x_n]: \sigma .f = f \,\,\, \forall \sigma \in S_{n}\}$.  Let $\Lambda(n)_k$ be the space of homogeneous symmetric polynomials of degree $k$. 
Then $\Lambda (n)$ is a graded ring with $\Lambda (n) = \bigoplus_k \Lambda (n)_k$. 

If we have a countable infinite set of indeterminates, say $x_1, x_2, \cdots $, we consider the ring $R$ of power
series in $x_1, x_2, \cdots $ of bounded degree. Hence, elements of $R$ can be infinite sums, but only in a finite number of degrees. Let $S_{\infty}$ be the group of all permutations of $\{1, 2, \cdots \}$. Then $S_{\infty}$ acts on $R$, and we define the ring of symmetric functions $\Lambda :=\{f \in R: \sigma . f=f \,\, for \,\, all \,\, \sigma \in S_{\infty}\}$. This is a subring of $R$. 
For each $n$ there is a surjective homomorphism $\Lambda \rightarrow \Lambda(n)$ obtained by setting $x_{n+1}=x_{n+2}= \cdots =0$. 
We also have a surjective homomorphism of rings $\Lambda (n+1) \rightarrow \Lambda (n)$ by setting $x_{n+1}=0$ and that restricts to a surjective map $\Lambda (n+1)_k \rightarrow \Lambda(n)_k$ for all $k$ and the map is bijective if and only if $k \leq n$.  We set $\displaystyle{\Lambda_k=\varprojlim_k   \Lambda(n)_k}$. Then $\Lambda=\displaystyle{\bigoplus_{k \geq 0} \Lambda_k}$. If $A$ is any commutative ring, we write $\Lambda_A:=\Lambda \otimes_{\mathbb Z} A$ for the ring of symmetric functions with coefficients in $A$.

There are various $\mathbb{Z}$-bases of the ring, some of which we shall review. 
They all are indexed by partitions. A partition $\lambda$ is a (finite or infinite) weakly decreasing sequence $\lambda=(\lambda_1, \lambda_2, \cdots )$ of non-negative integers with finitely many non-zero terms. Let $\mathcal P$ denote the set of all partitions. Given two partitions $\lambda$ and $\mu$ we say that $\lambda \geq \mu$ if and only if $|\lambda|=|\mu|$ and $\lambda_1+\lambda_2 \cdots \lambda_r \geq \mu_1 +\mu_2+ \cdots +\mu_r$ for all $r \geq 1$. For a partition $\lambda=(\lambda_1, \lambda_2, \cdots )$, the conjugate partition is defined by $\lambda'=(\lambda_1', \lambda_2', \cdots )$, where $\lambda_i'$ is the number of $j$'s such that $\lambda_j \geq i$.  

For a partition $\lambda=(\lambda_1, \lambda_2, \cdots )$ we define $x^{\lambda}:=x_1^{\lambda_1}x_2^{\lambda_2}\cdots $. The monomial symmetric function $m_{\lambda}$
is the sum of all distinct monomials obtainable from $x^{\lambda}$ by permutations of the $x$'s. In particular, when $\lambda=(1^r)=: (\underbrace{1,1,1,\cdots ,1}_{r}, 0,0, \cdots )$ we have $\displaystyle{m_{(1^r)}=e_r=\sum_{i_1 <i_2< \cdots <i_r}  x_{i_1}x_{i_2}\cdots x_{i_r}}$, the $r$-th elementary symmetric polynomial. The $e_r$'s are algebraically independent over $\mathbb Z$. At the other extreme when $\lambda=(r):=(r,0,0, \cdots )$ we have $m_{(r)}=p_r=\sum x_i^r$, the $r$-th power sum. It is clear that every $f \in \Lambda$  is uniquely expressible as a finite linear combination of the $m_{\lambda}$'s, so that $(m_{\lambda})_{\lambda \in \mathcal P}$ is a $\mathbb Z$-basis of $\Lambda$. For each $r \geq 0$, the $r$-th complete symmetric polynomial $h_r$, is the sum of all monomials of total degree $r$ in the variables $x_1, x_2, \cdots$ so that $h_r=\sum_{|\lambda|=r} m_{\lambda}$. The map $\omega: \Lambda \rightarrow \Lambda$ defined by $e_{\lambda} \mapsto h_{\lambda}$ is an involution of graded rings and hence $h_r$'s are also algebraically independent over $\mathbb Z$. 

For any partition $\lambda$ we define $e_{\lambda}:=e_{\lambda_1}e_{\lambda_2}\cdots $, $p_{\lambda}:=p_{\lambda_1}p_{\lambda_2} \cdots $ and $h_{\lambda}:=h_{\lambda_1}h_{\lambda_2} \cdots $ as the elementary symmetric function, the power sum symmetric function and the complete symmetric function respectively for the partition $\lambda$. It is easy to see that \[e_{\lambda'}=m_{\lambda}+\sum_{\mu < \lambda} a_{\lambda, \mu} m_{\mu}\] for some nonnegative integers $a_{\lambda, \mu}$. Hence $(e_{\lambda})_{\lambda \in \mathcal P}$ form another $\mathbb Z$-basis of $\Lambda$. Again since $\omega$ is an involution, the set $(h_{\lambda})_{\lambda \in \mathcal P}$ is yet another $\mathbb Z$-basis of $\Lambda$. 

The generating function for the $e_r$, $h_r$ and $p_r$ are $E(t)=\sum_{r \geq 0} e_rt^r=\prod (1+x_it)$, $H(t)=\sum_{r \geq 0} h_rt^r=\prod (1-x_it)^{-1}$ and 
$P(t)=\sum_{r \geq 1} p_rt^{r-1}=\sum_{r \geq 1} \frac{x_i}{(1-x_it)}$ respectively and they satisfy the following identities. \[E(t)H(-t)=1, P(t)=\frac{H'(t)}{H(t)} \,\, and \,\, P(-t)=\frac{E'(t)}{E(t)}.\] As a result for each $n \geq 1$ we get that 

\[\sum_{r=0}^n (-1)^r e_rh_{n-r}=0,\] \[nh_n=\sum_{r=1}^n p_rh_{n-r}\] and 

\begin{equation}{\label{e1}}
ne_n=\sum_{r=1}^n (-1)^{r-1}p_re_{n-r}.
\end{equation}

From the above expression it follows that $\mathbb Q[h_1, h_2, \cdots, h_n]=\mathbb Q[p_1, p_2,\cdots ,p_n]$ for all $n \geq 1$. Letting $n \rightarrow \infty$ we get that $\Lambda_{\mathbb Q}= \mathbb Q[h_1,h_2, \cdots ]=\mathbb Q[p_1, p_2,\cdots ]$. So $(p_{\lambda})_{\lambda \in \mathcal P}$ form a $\mathbb Q$-basis of $\Lambda_{\mathbb Q}$ but not a $\mathbb Z$-basis of $\Lambda$. 

We know that  $P(t) = \frac{d}{dt}\, log\, H(t)$. Integrating both sides imposing the boundary condition $log\, H(0) = 0$ and applying the exponential map, we get that, \[H(t)=\text{exp} \,\, (\sum_{n \geq 1} \frac{p_nt^n}{n})=\prod_{n \geq 1} \text{exp} \,\,  (\frac{p_nt^n}{n})=\prod_{n \geq 1}\sum_{d \geq 0} \frac{p_n^dt^{nd}}{n^dd!}=\sum_{\lambda} \frac{p_{\lambda}}{z_{\lambda}}t^{|\lambda|},\]
where $z_{\lambda}= \prod_{i \geq 1} i^{m_i(\lambda)}m_i(\lambda)!$ and $m_i(\lambda)$ is the number of $i$ appears in $\lambda$.

Since the involution $\omega$ maps $H(t)$ to $E(t)$ and vice versa, we have $\omega(p_n)=(-1)^{n+1}p_n$ and as a result, for a partition $\lambda$ we get that, $\omega(p_{\lambda})=\epsilon_{\lambda}p_{\lambda}$, where 
 $\epsilon_{\lambda}=(-1)^{|\lambda|-l(\lambda)}$. Now by applying the involution 
$\omega$ to the above identity we get that,

\[E(t)=\text{exp} \,\,(\sum_{n \geq 1} (-1)^{n+1}\frac{p_nt^n}{n})=\sum_{\lambda} \epsilon_{\lambda} \frac{p_{\lambda}}{z_{\lambda}}t^{|\lambda|}.\] 

Replacing $t$ by $-t$ in the last identity we get that 

\begin{equation}{\label{e2}}
E(-t)= \sum_{n \geq 0} (-1)^n e_nt^n= \text{exp} \,\, (-\sum_{n \geq 1} \frac{p_nt^n}{n})=\sum_{\lambda} (-1)^{(2|\lambda|-l(\lambda))} \frac{p_{\lambda}}{z_{\lambda}}t^{|\lambda|}.
\end{equation}

{\bf{Acknowledgments:}} Authors would like to thank S. Eswara Rao, who suggested this problem and for some helpful discussions. Authors would also like to thank S. Viswanath and Tanusree Khandai for some helpful discussions.

\end{document}